\documentclass{article}

\usepackage{amssymb}

\usepackage{graphicx}
\usepackage{amsmath}

\newfont{\bba}{eufm10 scaled\magstep{5}}
\newfont{\bbd}{eufm10 scaled\magstep2}
\newfont{\bbb}{msbm10 scaled \magstephalf}
\newfont{\bbc}{msbm8 scaled \magstephalf}

\begin{document}

\title{ SELECTIONS, PARACOMPACTNESS AND COMPACTNESS}

\author{Mitrofan M. Choban, Ekaterina P. Mihaylova, Stoyan I. Nedev}

\maketitle

\vskip 5pt

{\bf Abstract:} In the present paper, the Lindel\"of number and the
degree of compactness of spaces and of the cozero-dimensional kernel
of paracompact spaces are characterized in terms of selections of
lower semi-continuous closed-valued mappings into complete
metrizable (or discrete) spaces.

\vskip 5pt



\vskip 5pt

{\Large \bf Introduction}

\vskip 5pt

All considered spaces  are assumed to be $T_1$-spaces. Our
terminology comes, as a rule, from (\cite{Eng}, [7], \cite{Mich1},
\cite{RS}).

A topological space $X$ is called \emph{paracompact} if $X$ is
Hausdorff and every open cover of $X$ has a locally finite open
refinement.

One of the main results of the theory of continuous selections is
the following theorem:

{\bf Theorem 0.1 (E.Michael [9])} {\it For any lower semi-continuous
closed-valued mapping $\theta : X \rightarrow Y$ of a paracompact
space $X$ into
 a complete metrizable space $Y$ there exist a compact-valued
lower semi-continuous mapping $\varphi  :X \rightarrow Y$ and
 a compact-valued
upper semi-continuous mapping $\psi   :X \rightarrow Y$ such that
 $\varphi (x) \subseteq \psi (x) \subseteq \theta (x)$
for any $x \in X.$

Moreover, if $dim X = 0,$ then the selections $\phi , \psi $ of
$\theta $ are single-valued and continuous.}

It will be shown that the existence of upper semi-continuous
selections for lower semi-continuous closed-valued mappings into
 a discrete spaces implies the paracompactness of the domain (see [1 -5, 11]).

The aim of the present article is to the determine the conditions on
a space $X$ under which for any lower semi-continuous closed-valued
mapping $\theta : X \rightarrow Y$ of the space $X$ into
 a complete metrizable (or discrete) space $Y$ there exists a
selection $\varphi  :X \rightarrow Y$ for which the image $\varphi
(X)$ is "small" in a given sense.

A family $\gamma$ of subsets of a space $X$ is \emph{star-finite
(star-countable)} if for every element $\Gamma \in \gamma$ the set
$\{L \in \gamma: L \bigcap \Gamma \neq \emptyset \}$ is finite
(countable).

A topological space $X$ is called \emph{strongly paracompact or
hypocompact} if  $X$ is Hausdorff and every open cover of $X$ has a
star-finite open refinement.

The cardinal number $l(X) = min\{m:$ {\it every open cover of $X$
has an open refinement of cardinality} $\leq m\}$ is the
\emph{Lindel\"{o}f number }of $X$.

The cardinal number $k(X) = min \{ m : $ {\it every open cover of
$X$ has an open refinement of cardinality} $<m\}$ is \textit{the
degree of the compactness} of $X$.

Denote by $\tau^+ $ the least cardinal number larger than the
cardinal number  $\tau$. It is obvious that $l(X) \leq k(X) \leq
l(X)^+$.

For a space $X$ put $\omega (X) = \bigcup \{U:\ U$ {\it is open in
$X$ and} $dim U = 0\}$ and let $c \omega (X) = X \setminus  \omega
(X)$ be \textit{the cozero-dimensional kernel} of $X$ (See
\cite{ChMN}).

{\bf Lemma  0.2.} {\it Let $X$ be a paracompact space, $U$ be an
open subset of $X$ and $U\cap c\omega (X) \neq \emptyset$. Then $dim
\ cl_X (U\cap c \omega (X))\neq 0$.}

\textbf{Proof:} See \cite{ChMN}. $\Box$

A family $\xi $ of subsets of $X$ is called {\it $\tau $-centred} if
$\cap \eta \not=\emptyset $ provided $\eta \subseteq \xi $ and
$|\eta | < \tau .$

{\bf Lemma 0.3.} {\it Let $X$ be a paracompact space and $\tau $ be
an infinite cardinal. Then:

1. $l(X) \leq \tau $ if and only if any discrete closed subset of
$X$ has cardinality $\leq \tau .$

2. The following assertions are equivalent:

 a) $k(X) \leq \tau $.

 b) Any discrete closed subset of $X$ has cardinality $< \tau ;$

 c) $\cap \xi \not=\emptyset $ for any $\tau $-centered filter of closed subsets of $X$.
}

\textbf{Proof:} It is obvious. $\Box$

Assertions $2a$ and $2c$ are equivalent and implication $2a
\rightarrow 2b$ is true for every space $X$.

{\bf Lemma 0.4.} {\it Let $X$ be a metrizable space and $\tau $ be
an infinite not sequential cardinal. Then:

1. $l(X) \leq \tau $ if and only if  $w(X) \leq \tau .$

2. $k(X) \leq \tau $ if and only if $w(X)  < \tau .$ }

\textbf{Proof:} It is obvious. $\Box$

\vskip 5pt

\section{On the degree of compactness of spaces}

\vskip 5pt

A subset $L$ of a completely regular space $X$ is {\it bounded} in
$X$ if for every continuous function $f:X\rightarrow \mathbb{R}$ the
set $f(L)$ is bounded.

A space $X$ is called {\it $\mu $-complete} if it is completely
regular and the closure $cl_X L$ of every bounded subset $L$ of $X$
is compact.

Every paracompact space is $\mu$-complete. Moreover, every
Dieudonn\'{e} complete space is $\mu$-complete (see [6]).

{\bf Definition 1.1.} {\it Let $X$ be a space and $\tau$ be an
infinite cardinal. Put $k(X, \tau )=\bigcup \{U:U$ is open in $X$
and $k(cl_X U)<\tau\}$ and $c(X, \tau)= X\setminus k(X, \tau )$. For
every $x\in X$ put $k(x, X)= min\{k(cl_XU):U$ is an open in $X$
neighborhood of $x\}$.}

 By definition,  $k(X, \tau )= \{x\in X:k(x,X)<\tau\}$.

{\bf Lemma 1.2.}\label{Lemk(X)<tau} {\it Let $X$ be a space, $\tau$
be an infinite cardinal, $\{U_\alpha :\alpha \in A\}$ be an open
discrete family in $X$ and $k(X)\leq \tau $. Then:

1. $\mid A\mid <\tau$;

2. If $x_\alpha \in U_\alpha \cap k(X,\tau)$ for every $\alpha \in
A$, then $\sup \{k(x_\alpha, X):\alpha \in A\}<\tau$;

3. If $x_\alpha \in U_\alpha \cap c(X,\tau)$ for every $\alpha \in
A$, then $\mid A\mid < cf(\tau )$. }

\textbf{Proof:} Since $k(X)\leq \tau $, every discrete family in $X$
has cardinality $<\tau$.

Suppose that $x_\alpha \in U_\alpha \cap k(X,\tau)$ for every
$\alpha \in A$ and $\sup \{k(x_\alpha, X):\alpha \in A\}=\tau$. In
this case $\tau$ is a non-regular limit cardinal and $cf(\tau )\leq
\mid A\mid <\tau$. From our assumption it follows that there exists
a family of cardinals $\{\tau_\alpha:\alpha \in A\}$ such that
$\tau_\alpha < k(x_\alpha, X)$ for every $\alpha \in A$ and $\sup
\{\tau_\alpha :\alpha \in A\}=\tau$. For every $\alpha \in A$ there
exists an open family $\gamma_\alpha $ of $X$ such that $cl_X
U_\alpha \subseteq \bigcup \gamma_\alpha$ and $\mid \xi \mid \geq
\tau_\alpha $ provided $\xi \subseteq \gamma_\alpha$ and $cl_X
U_\alpha \subseteq \bigcup \xi$. One can assume that $U_\beta \cap
V=\emptyset$ for every $\alpha ,\ \beta \in A, \alpha \neq \beta$
and $V\in \gamma_\alpha$.
 Let $\gamma = (X\setminus \bigcup
\{cl_X U_\alpha :\alpha \in A\}) \cup (\bigcup  \{\gamma_\alpha
:\alpha \in A\} ).$  Then $\gamma$ is an open cover of $X$ and every
subcover of $\gamma $ has a cardinality $\geq
sup\{\tau_\alpha:\alpha \in A\}=\tau$, which is a contradiction.

If $x_\alpha \in U_\alpha \cap c(X,\tau)$ for every $\alpha \in A$
and $\mid A\mid \geq cf(\tau )$, then there exists a family of
cardinals $\{\tau_\alpha:\alpha \in A\}$ such that $\tau_\alpha <
\tau$ for every $\alpha \in A$ and $\sup \{\tau_\alpha :\alpha \in
A\}=\tau$. Since $k(x_\alpha, X)= \tau \geq \tau_\alpha$ for every
$\alpha \in A$, one can obtain a contradiction as in the previous
case. $\Box$

{\bf Lemma 1.3.}  {\it Let $X$ be a completely regular space, $\tau$
be a sequential cardinal and $k(X)\leq \tau$. Then the set $c(X,\tau
)$ is closed and bounded. Moreover, if $X$ is a $\mu$-complete
space, then:

1. $c(X,\tau)$ is a compact subset;

2. If $Y\subseteq k(X, \tau)$ is a closed subset of $X$, then $k(Y)<
\tau$. }

 \textbf{Proof:} If $\tau = \aleph_0$, then the space $X$ is compact
 and $k(X,\tau)$ is the subset of all isolated in $X$ points. Thus the
 set $c(X, \tau)$ is compact and every closed in $X$ subset of
 $k (X,\tau)$ is finite.

Suppose that $\tau$ is uncountable. There exists a family of
infinite cardinal numbers $\{\tau_n:n \in \mathbb{N}\}$ such that
$\tau_n < \tau_{n+1}<\tau$ for every $n \in \mathbb{N}$ and $\sup
\{\tau_n :n \in \mathbb{N}\}=\tau$. Suppose that the set $c(X,
\tau)$ is unbounded in $X$. Then there exists a continuous function
$f:X\rightarrow \mathbb{R}$ and a sequence $\{x_n\in c(X,\tau): n
\in \mathbb{N}\}$ such that $f(x_1)=1$ and $f(x_{n+1})\geq 3+f(x_n)$
for every $n \in \mathbb{N}$. The family $\xi = \{U_n=f^{-1}
((f(x_n)-1,f(x_n)+1)):n \in \mathbb{N}\}$ is discrete in $X$ and
$x_n\in U_n$ for every $n \in \mathbb{N}$. Then, by virtue of Lemma
1.2, $\mid \xi \mid < cf(\tau)=\aleph_0$, which is a contradiction.
Thus the set $c(X, \tau)$ is closed and bounded in $X$.

Assume now  that $X$ is a $\mu$-complete space. In this case the set
$c(X,\tau )$ is  compact.

Suppose that $Y\subseteq k(X, \tau)$ is a closed subset of $X$ and
$k(Y)= \tau$. We affirm that $\sup \{k(y, X):y\in Y\}<\tau$. For
every $x\in k(X,\tau)$ fix a neighborhood $U_x$ in $X$ such that
$k(cl_X U_x) =k(x,X)$. Suppose that $\sup \{k(y, X):y\in Y\}=\tau$.
For every $n\in \mathbb{N}$ fix a point $y_n \in Y$ such that
$k(y_n, X)\geq \tau_n$. Put $L=\{y_n:n\in \mathbb{N}\}$. If the set
$L$ is unbounded in $X$, then there exists a continuous function
$f:X\rightarrow \mathbb{R}$ such that $\sup \{f(y_n):n\in
\mathbb{N}\}=\infty$. One can assume that $f(y_{n+1})> 3+f(y_n)$.
The family $\xi = \{U_n=f^{-1} ((f(y_n)-1,f(y_n)+1)):n \in
\mathbb{N}\}$ is discrete in $X$ and $y_n\in U_n \cap k(X, \tau)$
for every $n \in \mathbb{N}$. Then, by virtue of Lemma 1.2, $\sup
\{k(y_n,X):n\in \mathbb{N}\}<\tau$, which is a contradiction. Thus
the set $L$ is bounded in $X$. Hence $cl_X L$ is a compact subset of
$Y$ and there exists an accumulation point $y\in cl_X L\setminus
(L\setminus \{y\})$. In this case $y\in k(X,\tau),\ k(y,X)<\tau$ and
$k(y,X)=\sup \{k(y_n,X):n\in \mathbb{N}\}=\tau$ which is a
contradiction. Hence $\sup \{k(y, X):y\in Y\}\leq \tau'<\tau$.

Since $k(X) \leq \tau $, then there exists a subset $Y'\subseteq Y$
such that $\mid Y'\mid \leq \tau '' $, $\tau'\leq \tau'' < \tau $
and $Y\subseteq \bigcup \{U_y:y\in Y'\}$. Since $k(cl_X U_y)\leq
\tau'\leq\tau''$ for every $y\in Y'$ and $\mid Y'\mid \leq \tau ''
$, then $l(\bigcup \{cl_X U_y:y\in Y'\})\leq \tau''<\tau$. Thus
$l(Y)\leq\tau''<\tau$ and $k(Y)<\tau$. $\Box$

 A subspace $Z$ of a space $X$ {\it is paracompact in} $X$ if for every
open family $\gamma = \{W_\mu : \mu \in M\}$ of $X$, for which
$Z\subseteq \cup \gamma ,$ there exists an open locally finite
family $\eta  = \{W_\mu' : \mu \in M\}$ of $X$ such that $Z\subseteq
\cup \eta  $ and $W_\mu ' \subseteq W_\mu $ for any $\mu \in M$.

{\bf Lemma 1.4.} \label{LemParacompact} {\it Let $X$ be a regular
space, $Z$ be a paracompact in $X$ subspace,
  $\tau$ be a limit cardinal number,
$k(Z)\leq \tau$ and $k(Y) < \tau $ for every closed subspace
$Y\subseteq X\setminus Z$ of the space $X$. Then:

1. $k(X) \leq \tau $, $c(X, \tau ) \subseteq Z$ and
$k(c(X,\tau))\leq cf(\tau )$;

2. If $Y\subseteq k(X, \tau)$ is a closed subset of $X$, then $k(Y)<
\tau$.

3. $Z$ is a closed subspace of $X$.}

\textbf{Proof:} Assertion 3 is obvious.

Let $x \in X\setminus Z$. Fix an open subset $U$ of $X$ such that $x
\in U \subseteq cl_XU \subseteq X\setminus Z.$ Then $k(cl_XU) < \tau
$ and  $c(X, \tau ) \subseteq Z$.

Let $\gamma $ be an open cover of $X$. Since $k(Z) \leq \tau $,
there exists a subsystem $\xi $ of $\gamma $ such that $|\xi |< \tau
$ and $Z \subseteq \cup \xi .$ Let $Y = X\setminus \cup \xi .$ Since
$k(Y) < \tau $, there exists a subsystem $\zeta  $ of $\gamma $ such
that $|\zeta  |< \tau $ and $Y \subseteq \cup \zeta  .$ Put $\eta
=\zeta \cup \xi .$ Then $\eta $ is a subcover of $\gamma $ and
$|\eta | < \tau .$ Thus $k(X) \leq \tau $.

 Suppose that $k(c(X,\tau))> cf(\tau )$. Since $Z$ is
paracompact in $X$, the subspace $c(X,\tau)$ is paracompact in $X$
and there exists an open locally-finite family $\{V_\alpha:\alpha
\in A\}$ of $X$ such that $c(X,\tau)\subseteq \bigcup \{V_\alpha :
\alpha \in A\},  |A| \geq cf(\tau)$ and $c(X,\tau)\setminus \bigcup
\{V_\beta : \beta \in A\setminus \{\alpha \}\}\neq \emptyset$ for
every $\alpha \in A$. For every $\alpha \in A$ fix $y_\alpha \in
c(X,\tau)\setminus \bigcup \{V_\beta : \beta \in A\setminus \{\alpha
\}\}\neq \emptyset$. Then $\{y_\alpha:\alpha \in A\}$ is a closed
discrete subset of $X$. There exists an open discrete family
$\{W_\alpha:\alpha \in A\}$ such that $y_\alpha \in W_\alpha
\subseteq V_\alpha$ for every $\alpha \in A$. By virtue of Lemma
1.2, $\mid A\mid < cf(\tau)$, which is a contradiction. Thus
$k(c(X,\tau))\leq cf(\tau )$.

Fix now a closed subset $Y$ of the space $X$ such that $Y\subseteq
k(X, \tau)$. We put $S = Y\cap Z$ and $\tau ' = \sup\{k(y,X):y\in
S\}$.

Suppose that $\tau ' = \tau .$ There exists a family of cardinals
$\{\tau_\alpha:\alpha \in A\}$ such that $\mid A\mid =cf(\tau)$,
$\sup \{\tau_\alpha :\alpha \in A\}=\tau$ and $\tau_\alpha <\tau$
for every $\alpha \in A$. One can assume that $A$ is well ordered
and $\tau_\alpha <\tau_\beta$ for every $\alpha, \beta \in A $ and
$\alpha <\beta$. For every $\alpha \in A$ there exists $y_\alpha \in
S$ such that $k(y_\alpha, X)>\tau_\alpha$. Let $L=\{y_\alpha :\alpha
\in A\}$. The cardinal $cf(\tau)$ is regular. If $y\in X$ and
$|W\cap L| = |A|$ for every neighborhood $W$ of $y$ in $X$, then
$y\in Y\subseteq k(X,\tau), k(y,X)<\tau$ and $k(y,X)\geq
\sup\{k(y_\alpha, X):\alpha \in A\}=\sup \{\tau_\alpha :\alpha \in
A\}=\tau$, which is a contradiction. Thus for every $y\in X$ there
exists an open neighborhood $W_y$ of $y$ in $X$ such that $|cl_X
W_y\cap L| < |A| = cf(\tau)$. There exists an open locally-finite
family $\{H_z:z \in Z\}$ of $X$ such that $Z \subseteq \cup  \{H_z:z
\in Z\}$ and $H_z\subseteq W_z$ for every $z\in Z$. Let $Z'=\{z\in
Z:H_z\cap L\neq \emptyset\}$. The set $Z'$ is discrete and closed in
$X$. Since $Z$ is a paracompact space, we have $|Z'| = \tau '' <
\tau $ and $|H_z\cap L| = \tau (z) < cf(\tau )$ for any $z \in Z'$.
Thus $|L| = \mid \cup \{H_z\cap L: z \in Z'\}\mid < \tau $, a
contradiction. Therefore $\tau ' < \tau $.

Since $S$ is paracompact in $X$ and $S\subseteq k(X,\tau )$,
 there exist a set $M \subseteq S$ and
a locally finite open in $X$ family $\{U_\mu : \mu \in M\}$ such
that $k(cl_X U_\mu ) \leq \tau '$, $|M| = \tau _1 <\tau $ and $S
\subseteq \cup \{U_\mu :\mu \in M\}$. Then $k(S_1) = \tau _2 <
\tau$, where $S_1 = \cup \{cl_X U_\mu : \mu \in M\}$. Let $Y_1 =
Y\setminus \cup \{U_\mu :\mu \in M\}$. Since $Y_1$ is a closed
subset of $X$ and $Y_1 \subseteq X\setminus Z$, $k(Y_1) = \tau _3 <
\tau $. Thus $k(Y) \leq k(Y_1) + k (S_1) <\tau $. $\Box$

{\bf Corollary 1.5.} \label{LemParacompact} {\it Let $X$ be a
paracompact space, $\tau$ be a limit cardinal and $k(X)\leq \tau$.
Then:

1. $k(c(X,\tau))\leq cf(\tau )$;

2. If $Y\subseteq k(X, \tau)$ is a closed subset of $X$, then $k(Y)<
\tau$.}

A {\it shrinking} of a cover $\xi = \{U_{\alpha } :\  \alpha \in
A\}$ of the space $X$ is a cover $\gamma = \{V_{\alpha } :\  \alpha
\in A\}$ such that $V_{\alpha } \subseteq U_{\alpha }$ for every
$\alpha \in A$ (see \cite{Eng}, \cite{Eng2}). The operation of
shrinking preserves the properties of local finiteness,
star-finiteness and star-countableness.

Let $\tau$ be an infinite cardinal number. A family $\gamma$ of
subsets of a space $X$ is called $\tau$-\textit{star}
($\tau^-$-\textit{star}) if $\mid \{H\in \gamma :\ H\cap L\neq
\emptyset \}\mid \leq \tau $ ($\mid \{H\in \gamma :\ H\cap L\neq
\emptyset \}\mid < \tau $) for every $L\in \gamma$.

A family $\{H_{\alpha } : \alpha \in A\}$ of subsets of a space $X$
is \textit{closure-preserving} if  $\bigcup \{cl_X H_{\beta } :
\beta \in B\} = cl_X ( \bigcup \{ H_{\beta } : \beta \in B\})$ for
every $B \subseteq A$ (see \cite{Mich3}).

{\bf Proposition 1.6.} {\it Let $\tau$ be an infinite cardinal and
$X$ be a paracompact space. Then the following assertions are
equivalent:

1. $k(c\omega (X)) \leq \tau$.

2. For every open cover of $X$ there exists an open $\tau^-$-star
shrinking.

3. For every open cover of $X$ there exists a closed
closure-preserving $\tau^-$-star shrinking.

4. For every open cover of $X$ there exists a closed  $\tau^-$-star
shrinking.}

\textbf{Proof:}  $(1\Rightarrow 2)$ and $(1\Rightarrow 3)$  Let
$\xi=\{U_\alpha:\ \alpha \in A\}$ be an open cover of $X$. There
exist a subset $B$ of $A$ and an open-and-closed subset $H$ of $X$
such that $c\omega (X) \subseteq H \subseteq \bigcup \{U_\alpha :\
\alpha \in B\}$ and $\mid B\mid <\tau $ (see the proof of
Proposition 4 \cite{ChMN}). Since $dim (X\setminus H)=0$ (unless
$X\setminus H$ is empty) there exists a discrete family $\{W_\alpha
:\ \alpha \in A\}$ of open-and-closed subsets of $X$ such that
$\bigcup \{W_\alpha :\ \alpha \in A\}=X\setminus H$ and $W_\alpha
\subseteq U_\alpha$ for every $\alpha \in A$. Let $V_\alpha =
(U_\alpha \cap H)\cup W_\alpha$ for $\alpha \in B$ and $V_\alpha =
W_\alpha $ for $\alpha \in A\setminus B$. Obviously $\gamma
=\{V_\alpha:\ \alpha \in A\}$ is an open $\tau^-$-star shrinking of
$\xi$.

Since $X$ is paracompact, there exists a closed locally finite
family $\{H_\alpha : \alpha \in B\}$ such that $H = \cup \{H_\alpha
: \alpha \in B\}$ and $H_\alpha \subseteq U_\alpha $ for any $\alpha
\in B$.  Put $H_\alpha = W_\alpha $ for any $\alpha  \in A\setminus
B.$
 Obviously $\lambda  =\{H_\alpha:\ \alpha
\in A\}$ is a closed locally finite $\tau^-$-star shrinking of
$\xi$. Every locally finite family is closure-preserving.
Implications $(1\Rightarrow 2)$ and $(1\Rightarrow 3)$  are proved.

Implication $(3\Rightarrow 4)$ is obvious.

$(2\Rightarrow 1)$ and $(4\Rightarrow 1)$ Suppose $k(c\omega (X))>
\tau$. There exists a locally finite open cover $\xi=\{U_\alpha:\
\alpha \in A\}$ of $c\omega (X)$ such that $c\omega (X) \setminus
\bigcup \{U_\alpha:\ \alpha \in B\}\neq \emptyset$ provided
$B\subseteq A$ and $\mid B\mid <\tau$. One can assume that $c\omega
(X) \setminus \bigcup \{U_\alpha:\ \alpha \in B\}\neq \emptyset$ for
every proper subset $B$ of $A$. Fix a point $x_\alpha \in c\omega
(X) \setminus \bigcup \{U_\beta:\ \beta \in A\setminus \{\alpha
\}\}$ for every $\alpha \in A$. The set $\{x_\alpha:\ \alpha \in
A\}$ is discrete in $X$. There exists a discrete family $\{V_\alpha
:\ \alpha \in A\}$ of open subsets of $X$ such that $x_\alpha \in
V_\alpha \subseteq cl_X V_\alpha \subseteq U_\alpha$ for every
$\alpha \in A$. Let $X_\alpha = cl_X V_\alpha $. Then $dim X_\alpha
>0$ and there exist two closed disjoint subsets $F_\alpha$ and
$P_\alpha $ of $X_\alpha$ such that if $W_\alpha $ and $O_\alpha$
are open in $X$ and $F_\alpha \subseteq W_\alpha \subseteq
X\setminus P_\alpha $, $P_\alpha \subseteq O_\alpha \subseteq
X\setminus F_\alpha$ and $X_\alpha \subseteq W_\alpha \cup O_\alpha
$, then
 $X_\alpha \cap W_\alpha \cap O_\alpha \not=\emptyset$.
The family $\{F_\alpha:\ \alpha \in A\}$ and the family
$\{P_\alpha:\ \alpha \in A\}$ are discrete in $X$. There exists a
discrete family $\{Q_\alpha:\ \alpha \in A\}$ of open subsets of $X$
such that $(\bigcup \{Q_\alpha :\ \alpha \in A\}) \cap (\bigcup
\{F_\alpha :\ \alpha \in A\})=\emptyset$, $P_\alpha \subseteq
Q_\alpha$ and $Q_\alpha \cap (\bigcup \{X_\beta :\ \beta \in
A\setminus \{\alpha \}\})=\emptyset$ for every $\alpha \in A$. Let
$\mu \notin A$, $M=A\cup\{\mu\}$ and $Q_\mu = X\setminus \bigcup
\{P_\alpha :\ \alpha \in A\})$. Then $\zeta =\{Q_m:\ m \in M\}$ is
an open cover of $X$. If $\gamma=\{H_m:\ m\in M\}$ is an open
shrinking of $\zeta$, then $H_\mu \cap H_\alpha \neq \emptyset$ for
every $\alpha \in A$. The last contradicts 2. Suppose now that
$\gamma=\{H_m:\ m\in M\}$ is a closed shrinking of $\zeta$. Let
$\alpha \in A$ and $H_\alpha \cap H_\mu = \emptyset .$ There exist
two disjoint open subsets $W_\alpha $ and $O_\alpha $ of $X$ such
that $H_\alpha \subseteq W_\alpha $ and $H_\mu \subseteq O_\alpha .$
Then $X_\alpha \subseteq H_\alpha \cup H_\mu \subseteq O_\alpha \cup
W_\alpha $, $P_\alpha \subseteq W_\alpha \subseteq X\setminus
F_\alpha $, $F_\alpha \subseteq O_\alpha \subseteq X\setminus
P_\alpha$, $X_\alpha \subseteq W_\alpha \cup O_\alpha $ and
 $X_\alpha \cap W_\alpha \cap O_\alpha =\emptyset$.
The last contradicts 4. Implications  $(2\Rightarrow 1)$ and
$(4\Rightarrow 1)$ are proved.$\Box$

\vskip 5pt

\section{The degree of compactness and selections}

\vskip 5pt

Let $X$ and $Y$ be non-empty topological spaces. A
\textit{set-valued mapping} $\theta:X\rightarrow Y$  assigns to
every $x \in X$ a non-empty subset $\theta(x)$ of $Y$. If $\phi, \
\psi: X\rightarrow Y$ are set-valued mappings and $\phi (x)
\subseteq \psi (x)$ for every $x \in X$, then $\phi$ is called
\textit{a selection} of $\psi$.

Let $\theta :X\rightarrow Y$ be a set-valued mapping and let $A
\subseteq X$ and $B\subseteq Y$. The set $\theta^{-1}(B)=\{x\in X:
\theta (x) \bigcap B \neq\emptyset \}$ is \textit{the inverse image}
of the set $B$, $\theta (A) =\theta^1 (A) = \bigcup \{\theta (x):x
\in A\}$ is \textit{the image} of the set $A$ and
$\theta^{n+1}(A)=\theta (\theta^{-1}(\theta^n (A))) $ is
\textit{the} $n+1$-\textit{image} of the set $A$. The set
$\theta^{\infty} (A)=\bigcup \{\theta^n (A):n \in \mathbb{N}\}$ is
\textit{the largest image} of the set A.

A set-valued mapping  $\theta :X\rightarrow Y$  is called
\textit{lower (upper) semi-continuous} if for every open (closed)
subset $H$ of $Y$ the set $\theta^{-1} (H)$ is open (closed) in $X$.

In the present section we study the mutual relations between the
following properties of topological spaces:

$K1.$ $k(X) \leq \tau$.

$K2.$ {\it For every lower semi-continuous closed-valued mapping
$\theta :X \rightarrow Y$ into a complete metrizable space $Y$ there
exists a lower semi-continuous selection $\phi :X \rightarrow Y$ of
$\theta $ such that $k(cl_{Y } \phi (X)) \leq \tau $.

$K3.$ For every lower semi-continuous closed-valued mapping $\theta
:X \rightarrow Y$ into a complete metrizable space $Y$ there exists
a set-valued selection $g :X \rightarrow Y$ of $\theta $ such that
$k(cl_{Y } g (X)) \leq \tau $.

$K4.$ For every lower semi-continuous closed-valued mapping $\theta
:X \rightarrow Y$ into a complete metrizable space $Y$ there exists
a single-valued selection $g :X \rightarrow Y$ of  $\theta $ such
that $k(cl_{Y } g (X)) \leq \tau $.

$K5.$  For every lower semi-continuous  mapping $\theta :X
\rightarrow Y$ into a discrete  space $Y$ there exists a lower
semi-continuous selection $\phi :X \rightarrow Y$ of  $\theta $ such
that $ |\phi (X)| < \tau $.

$K6.$ For every lower semi-continuous mapping $\theta :X \rightarrow
Y$ into a discrete space $Y$ there exists a set-valued selection $g:
X \rightarrow Y$ of  $\theta $ such that $\mid g (X)\mid < \tau $.

$K7.$ For every lower semi-continuous mapping $\theta :X \rightarrow
Y$ into a discrete space $Y$ there exists a single-valued selection
$g:X \rightarrow Y$ of  $\theta $ such that $\mid g (X)\mid < \tau
$.

$K8.$ Every open cover of $X$ has a subcover of cardinality $< \tau.
$

$K9$. For every lower semi-continuous closed-valued mapping $\theta
:X \rightarrow Y$ into a complete metrizable space $Y$ there exist a
compact-valued lower semi-continuous mapping $\varphi  :X
\rightarrow Y$ and
 a compact-valued
upper semi-continuous mapping $\psi   :X \rightarrow Y$ such that
$k(cl_Y(\psi (X))) \leq \tau $ and $\varphi (x) \subseteq \psi (x)
\subseteq \theta (x)$ for any $x \in X.$

$K10.$ \it For every lower semi-continuous closed-valued mapping
$\theta :X \rightarrow Y$ into a complete metrizable space $Y$ there
exists an upper semi-continuous selection $\phi :X \rightarrow Y$ of
$\theta $ such that $k(cl_{Y } \phi (X)) \leq \tau $.

$K11.$ \it For every lower semi-continuous closed-valued mapping
$\theta :X \rightarrow Y$ into a complete metrizable space $Y$ there
exists a lower semi-continuous selection $\phi :X \rightarrow Y$ of
$\theta $ such that $ w(\phi (X)) < \tau $.

$K12.$  \it For every lower semi-continuous closed-valued mapping
$\theta :X \rightarrow Y$ into a complete metrizable space $Y$ there
exist a closed-valued lower semi-continuous selection $\phi  :X
\rightarrow Y$ of  $\theta $, a selection $\mu   :X \rightarrow Y$
of  $\theta $ and a closed $G_\delta $ subset $F$ of the space $X$
such that:

- $\phi (x) \subseteq \mu (x)$ for any $x \in X$;

- $\phi (x) = \mu (x)$ for any $x \in X\setminus F$;

- the mapping $\mu |F  :F \rightarrow Y$ is upper semi-continuous
and is closed-valued;

- $c(X,\tau ) \subseteq F$, $\Phi = \mu (F)$ is a compact subset of
$Y$ and $k(Z\cap \mu (X)) < \tau $ provided $Z \subseteq Y\setminus
\Phi $ and $Z$ is a closed subspace of the space $Y$;

- $k(cl_Y \phi (X)) \leq k(\mu (X)) \leq \tau $.}

Let us mention that, in the conditions of $K12$:

 - $k(\mu (X)) \leq \tau $ provided the set $\Phi $ is compact and $k(Z\cap \mu (X)) < \tau $
 for a closed subset $Z \subseteq Y\setminus \Phi $ of the space $Y$;

 - the mapping $\mu   :X \rightarrow Y$ is closed-valued and the
mapping $\mu |F  :F \rightarrow Y$ is compact-valued;

- the mapping $\mu |(X\setminus F)  :X\setminus F \rightarrow Y$ is
lower semi-continuous;

- the mapping $\mu   :X \rightarrow Y$ is Borel measurable, i.e.
$\mu ^{-1}(H)$ is a Borel subset of the space $X$ for any open or
closed subset $H$ of $Y$.

The $\sigma $-algebra generated by the open subsets of the space $X$
is the algebra of Borel subsets of the space $X$.

{\bf Lemma 2.1.} {\it  Let $X$ be a space and $\tau $ be an infinite
cardinal. Then the following implications $(K9 \rightarrow K2
\rightarrow K3$ $ \rightarrow K4 \rightarrow K3 $ $\rightarrow K6
\rightarrow  K7 \rightarrow K8 \rightarrow K1 \rightarrow K5
\rightarrow K6$, $K10 \rightarrow K3)$ and $(K12 \rightarrow K11
\rightarrow  K2 \rightarrow K5 \rightarrow K6)$ are true.}

\textbf{Proof:} Implications $(K12 \rightarrow K11 \rightarrow K2
\rightarrow K5 \rightarrow K6$, $K9 \rightarrow K2 \rightarrow K3
\rightarrow K6)$, $(K4 \rightarrow K7$, $K4\rightarrow K3)$, $(K7
\rightarrow K6$, $K8 \rightarrow K1 \rightarrow K8)$ and $(K10
\rightarrow  K3)$ are obvious.

Let $\phi:X \rightarrow Y$ be a set-valued selection of the mapping
$\theta :X \rightarrow Y$ and $k(cl_Y \phi (X))\leq \tau$. For every
$x\in X$ fix a point $f(x)\in \phi (x)$. Then $f:X \rightarrow Y$ is
a single-valued selection of $\theta $ and $\phi$, $f(X)\subseteq
\phi(X)$ and $k(cl_Y f (X))\leq  k(cl_Y \phi (X))\leq \tau$. The
implications $(K3 \rightarrow  K4)$ and $(K6 \rightarrow  K7)$ are
proved.

Let $\gamma =\{U_\alpha:\alpha \in A\}$ be an open cover of $X$. One
may assume that $A$ is a discrete space. For every $x\in X$ put
$\theta_\gamma (x) =\{\alpha \in A:x\in U_\alpha\}$. Since
$\theta^{-1} (\{\alpha \})=U_\alpha$, the mapping $\theta_\gamma $
is lower semi-continuous. Let $\phi :X \rightarrow Y$ be a
set-valued selection of $\theta_\gamma $ and $\mid \phi (X)
\mid<\tau$. Put $B=\phi (X)$ and $H_\alpha=\phi^{-1} (\alpha )$ for
every $\alpha \in B$. Then $H_\alpha \subseteq \theta^{-1} (\{\alpha
\})=U_\alpha$ for every $\alpha \in B$, $X=\bigcup \{H_\alpha:\alpha
\in B\}$, $\xi=\{H_\alpha:\alpha \in B\}$ is a refinement of $\gamma
$ and $\mid B\mid <\tau$. Implications $(K3 \rightarrow  K8)$ and
$(K6 \rightarrow  K8)$ are proved.

Let $k(X)\leq \tau$ and $\theta :X \rightarrow Y$ be a lower
semi-continuous mapping into a discrete space $Y$. Then
$\{U_y=\theta^{-1} (y):y\in Y\}$ is an open cover of $X$. There
exists a subset $Z\subseteq Y$ such that $\mid Z\mid <\tau$ and
$X=\bigcup \{U_y:y\in Z\}$. Now we put $\phi (x)=\{y\in Z:x \in
U_y\}$. Then $\phi :X \rightarrow Y$ is a lower semi-continuous
selection of $\theta $, $\phi (x)=Z\cap \theta (x)$ for every $x\in
X$ and $\mid \phi(X)\mid =\mid Z\mid <\tau$. Implication $(K1
\rightarrow  K5)$ is proved. The proof is complete.

{\bf Proposition 2.2.} {\it Let $X$ be a space, $\tau$ be an
infinite cardinal and $\theta :X \rightarrow Y$ be an upper
semi-continuous mapping onto $Y$. Then:

1. If $l(X)\le \tau$ and $l(\theta (x))\leq \tau$ for every $x\in
X$, then $l(Y)\le \tau$;

2. If $k(X)\le \tau$ and $k(\theta (x))\leq cf(\tau)$ for every
$x\in X$, then $k(Y)\le \tau$;

3. If $\theta $ is compact-valued, then $l(Y)\leq l(X)$ and
$k(Y)\leq k(X)$.

4. If $X$ is a $\mu $-complete space, $\tau $ is a sequential
cardinal number and $\theta $ is compact-valued, then $c(Y,\tau )
\subseteq \theta (c(X,\tau  ))$ and $k(Z) < \tau $ provided
$Z\subseteq Y\setminus c(Y,\tau )$ and $Z$ is closed in the space
$Y$.}

\textbf{Proof:} If $V$ is an open subset of $Y$, then $\theta^*
(V)=\{x\in X:\theta (x)\subseteq V\}$ is open in $X$.

1. Let $\tau$ be an infinite cardinal,  $l(X)\le \tau$ and $l(\theta
(x))\leq \tau$ for every $x\in X$. Let $\gamma =\{V_\alpha :\alpha
\in A\}$ be an open cover of $Y$. If $x\in X$, then $l(\theta
(x))\leq \tau$. Thus every open family in $Y$, which covers $\theta
(x)$, has a subfamily of cardinality $\leq \tau$ covering $\theta
(x)$. Hence there exists a subset $A_x\subseteq A$ such that $\mid
A_x \mid = \tau_x \leq \tau$ and $\theta (x) \subseteq \bigcup
\{V_\alpha: \alpha \in A_x\}$. We put $W_x = \cup \{V_\alpha :
\alpha \in A_x\}$ and $U_x = \{z \in X: \theta (z) \subseteq W_x\}$.

Obviously $\lambda =\{U_x: x\in X\}$ is an open cover of $X$. Since
$l(X)\leq \tau$, there exists an open subcover $\zeta= \{U_x: x \in
X'\}$ of $\lambda $ such that $\mid X'\mid \leq \tau$ and
$X'\subseteq  X.$ Let $B = \cup \{A_x: x \in X'\}$. Obviously $|B|
\leq \tau .$ Since $\theta (U_x) \subseteq W_x$ for any $x \in X,$
we have $Y = \theta (X)= \theta (\cup \{U_x: x \in X'\})$ = $\cup
\{\theta (U_x) : x \in X'\} \subseteq \cup \{W_x: x \in X'\}$ =
$\cup \{V_\alpha : \alpha \in B\}$. Hence $\gamma ' = \{V_\alpha :
\alpha \in B\}$ is a subcover of $\gamma $ of cardinality $\leq \tau
.$ Assertion 1 is proved.

2. One can follow the proof of the previous assertion 1. Let $\tau$
be an infinite cardinal,
 $k(X)\le \tau$ and
$k(\theta(x))\leq cf(\tau)$ for every $x\in X$. Let $\gamma
=\{V_\alpha :\alpha \in A\}$ be an open cover of $Y$. For any $x \in
X$ there exists a subset $A_x\subseteq A$ such that $\mid A_x \mid =
\tau_x < cf(\tau )$ and $\theta (x) \subseteq \bigcup \{V_\alpha:
\alpha \in A_x\}$. We put $W_x = \cup \{V_\alpha : \alpha \in A_x\}$
and $U_x = \{z \in X: \theta (z) \subseteq W_x\}$.

Obviously $\lambda =\{U_x: x\in X\}$ is an open cover of $X$. Since
$k(X)\leq \tau$, there exists an open subcover $\zeta= \{U_x: x \in
X'\}$ of $\lambda $ such that $\mid X'\mid = \tau _0 < \tau$ and
$X'\subseteq  X.$ Let $B = \cup \{A_x: x \in X'\}$. Since $\theta
(U_x) \subseteq W_x$ for any $x \in X,$ we have $Y = \theta (X)=
\theta (\cup \{U_x: x \in X'\}) = \cup \{\theta (U_x) : x \in X'\}$
$\subseteq \cup \{W_x: x \in X'\}= \cup \{V_\alpha : \alpha \in
B\}$. Hence $\gamma ' = \{V_\alpha : \alpha \in B\}$ is a subcover
of $\gamma $.

We affirm that $|B| < \tau .$

Consider the following cases:

Case 1. $\tau$ is regular, i.e. $cf(\tau )=\tau$.

Since $\mid X'\mid =\tau_0 <\tau =cf(\tau )$ and $\mid A_x \mid  <
\tau $ for every $x\in X$, it follows that $|B| \leq \Sigma
\{\tau_{x}: x \in X'\} = \tau'<\tau $.

 Hence $\gamma ' = \{V_\alpha : \alpha \in B\}$ has cardinality $<\tau $.

Case 2. $\tau$ is not regular, i.e. $cf(\tau )= m <\tau$.

In this case $\tau $ is a limit cardinal, $\tau _0 <\tau $ and $m <
\tau $. Hence $\tau ' = sup\{m, \tau _0\} < \tau .$

Since $\mid A_x \mid = \tau_x < m$ for every $x\in X$, it follows
that $|B| \leq  \Sigma \{\tau_{x}: x \in X'\} \leq \tau ' < \tau $.

Hence $\gamma ' = \{V_\alpha : \alpha \in B\}$ has cardinality
$<\tau $.

Assertion 2 is proved.

3.  Assertion 3 follows easily from assertions 1 and 2.

4.Obviously, $\Phi = \theta (c(X,\tau ))$ and $c(Y,\tau )$ are
compact subsets of the space $Y$. Let $Z \subseteq Y\setminus \Phi $
be a closed subspace of the space $Y$. Then $X_1 = \theta ^{-1}(Z)$
is a closed subspace of the space $X$ and $X_1\cap c(X,\tau ) =
\emptyset $. By virtue of Lemma 1.3, $k(X_1)< \tau $. Let $Y_1 =
\theta (X_1)$. Then $\theta _1 = \theta |X_1: X_1\rightarrow Y_1$ is
 an upper semi-continuous mapping onto $Y_1$. From assertion 2
it follows that $k(Y_1) \leq k(X_1) < \tau $. Since $Z$ is a closed
subspace of the space $Y_1$, we have $k(Z) \leq k(Y_1) < \tau $. In
particular, $Y\setminus \Phi \subseteq k(Y,\tau )$ and $c(Y,\tau
)\subseteq \Phi $. Since $\Phi $ is a compact subset of $Y$, $k(Z)
<\tau $ provided $Z\subseteq Y\setminus c(Y,\tau )$ and $Z$ is
closed in the space $Y$.
 $\Box$

{\bf Theorem 2.3.}\label{NonSeq} {\it Let $X$ be a regular space and
$\tau $ be a regular  cardinal number. Then assertions $K1 - K8$ and
$K12$ are equivalent. Moreover, if the cardinal number $\tau $ is
regular and uncountable, then assertions $K1 - K8$, $K11$ and $K12$
are equivalent.}

\textbf{Proof:}  Let $k(X) \leq \tau$ and  $\theta : X \rightarrow
Y$ be a lower semi-continuous closed-valued mapping into a complete
metric space $(Y, \rho )$.

Case 1. $\tau =\aleph_0$.

In this case the space $X$ is compact. Thus, from E.Michael's
Theorem [9] (see Theorem 0.1), it follows that there exist a lower
semi-continuous compact-valued mapping $\varphi  :X \rightarrow Y$
and an upper semi-continuous compact-valued mapping $\psi  :X
\rightarrow Y$ such that $\varphi (x) \subseteq \psi (x) \subseteq
\theta (x)$ for any $x \in X.$ The set $\psi (X)$ is compact and
$\varphi (X) \subseteq \psi (X)$. Implication $(K1\Rightarrow K9)$
is proved.

Case 2. $\tau > \aleph_0$.

There exists a sequence $\gamma = \{\gamma _n =\{U_\alpha : \alpha
\in A_n\}: n \in \mathbb N\}$ of open covers of the space $X$, a
sequence $\xi  = \{\xi  _n =\{V_\alpha : \alpha \in A_n\}: n \in
\mathbb N\}$ of open families of the space $Y$ and a sequence $\pi =
\{\pi _n: A_{n+1} \rightarrow A_n: n \in \mathbb N\}$ of mappings
such that:

- $\cup \{U_\beta : \beta \in \pi _{n}^{-1}(\alpha )\}$ = $U_{\alpha
} \subseteq cl_X U_\alpha \subseteq \theta ^{-1}(V_\alpha )$ for any
$\alpha \in A_n$ and $n \in \mathbb N$;

- $\cup \{cl_{Y}V_\beta : \beta \in \pi _{n}^{-1}(\alpha )\}$
$\subseteq V_\alpha  $ and $diam(V_\alpha ) < 2^{-n}$ for any
$\alpha \in A_n$ and $n \in \mathbb N$;

- $|A_n| < \tau $ for any $n \in \mathbb N$.

Let $\eta = \{V: V$ {\it is open in $Y$ and} $diam (V) < 2^{-1}\}$.
Let $\gamma ' = \{U: U$ {\it is open in $X$ and $cl_{X}U \subseteq
\theta ^{-1}(V)$ for some} $V \in \eta  \}$. Since $k(X) \leq \tau
$, there exists an open subcover $\gamma _1 =\{U_\alpha : \alpha \in
A_1\}$ of $\gamma '$ such that $|A_1| < \tau .$ For any $\alpha \in
A_1$ fix $V_\alpha \in \eta $ such that $cl_X U_\alpha \subseteq
\theta ^{-1}(V_\alpha )$.

Consider that the objects $\{\gamma _i, \xi _i, \pi _{i-1}: i \leq
n\}$ are constructed. Fix $\alpha \in A_n$. Let $\eta_{\alpha } =
\{V: V$ {\it is open in $Y$, $cl_YV \subseteq V_\alpha $
 and} $diam (V) < 2^{-n-1}\}$.
Let $\gamma_{\alpha }' = \{W: W$ {\it is open in $X$ and $cl_{X}W
\subseteq \theta ^{-1}(V)$ for some} $V \in \eta_\alpha   \}$. Since
$k(cl_XU_\alpha ) \leq \tau $ and $cl_XU_\alpha \subseteq \cup
\gamma_{\alpha }'$, there exists an open subfamily $\gamma _\alpha
=\{W_\beta  : \beta  \in A_\alpha \}$ of $\gamma_\alpha  '$ such
that $|A_\alpha | < \tau $ and $cl_XU_\alpha \subseteq \cup
\{W_\beta :\beta \in A_\alpha \}.$ For any $\beta  \in A_\alpha $
fix $V_\beta  \in \eta _\alpha $ such that $cl_X W_\beta  \subseteq
\theta ^{-1}(V_\beta  )$. Let $A_{n+1} = \cup \{A_\alpha : \alpha
\in A_n\}$, $\pi _{n}^{-1}(\alpha ) = A_\alpha $ and $U_\beta =
U_\alpha \cap W_\beta $ for all $\alpha \in A_n$ and $\beta \in
A_\alpha .$ Since $\tau $ is regular and uncountable, then
$|A_{n+1}| < \tau .$

The objects $\{\gamma _n, \xi _n, \pi _{n}: n \in \mathbb N\}$ are
constructed.

Let $x \in X.$ Denote by $A(x)$ the set of all sequences $\alpha =
(\alpha _n: n \in \mathbb N)$ for which $\alpha _n \in A_n$ and $x
\in U_{\alpha _n}$ for any  $n \in \mathbb N.$ For any  $\alpha =
(\alpha _n: n \in \mathbb N) \in A(x)$ there exists a unique point
$y(\alpha ) \in Y$ such that $\{y(\alpha )\} = \cap \{V_{\alpha _n}:
n \in \mathbb N\}.$ It is obvious that $y(\alpha ) \in \theta (x).$
Let $\phi (x) = \{y(\alpha ): \alpha \in A(x)\}.$ Then $\phi $ is a
selection of $\theta $. By construction:

- $U_\alpha \subseteq \phi ^{-1}(V_\alpha )$ for all $\alpha \in
A_n$ and $n \in \mathbb N$;

- the mapping $\phi $ is lower semi-continuous;

- if $Z= \phi (X)$, then $\{H_\alpha = Z\cap V_\alpha : \alpha \in A
= \cup \{A_n:  n \in \mathbb N\}\}$ is an open base of the subspace
$Z$.

We affirm that $ w(Z) < \tau $.

Subcase 2.1. $\tau $ is a limit cardinal.

In this subcase $ m = sup \{|A_n|: n \in \mathbb N\} < \tau $ and
$w(Z) \leq |A| \leq m <\tau $.

Subcase 2.2. $\tau $ is not a limit cardinal.

In this subcase there exists a cardinal number $m$ such that $m ^+ =
\tau $ and $|A| \leq m$. Thus $w(Z) < \tau $.

In this case we have proved implication $(K1 \rightarrow K11).$

Lemma 2.1 completes the proof of the theorem. $\Box$

{\bf Corollary 2.4.} {\it Let $X$ be a regular space and $\tau $ be
a  cardinal number. Then the following assertions  are equivalent:

$L1.$ $l(X) \leq \tau$.

$L2.$  For every lower semi-continuous closed-valued mapping $\theta
:X \rightarrow Y$ into a complete metrizable space $Y$ there exists
a lower semi-continuous selection $\phi :X \rightarrow Y$ of
$\theta $ such that $l(cl_{Y } \phi (X)) \leq \tau $.

$L3.$ For every lower semi-continuous closed-valued mapping $\theta
:X \rightarrow Y$ into a complete metrizable space $Y$ there exists
a set-valued selection $g :X \rightarrow Y$ of $\theta $ such that
$l(cl_{Y } g (X)) \leq \tau $.

$L4.$ For every lower semi-continuous closed-valued mapping $\theta
:X \rightarrow Y$ into a complete metrizable space $Y$ there exists
a single-valued selection $g :X \rightarrow Y$ of  $\theta $ such
that $l(cl_{Y } g (X)) \leq \tau $.

$L5.$  For every lower semi-continuous  mapping $\theta :X
\rightarrow Y$ into a discrete  space $Y$ there exists a lower
semi-continuous selection $\phi :X \rightarrow Y$ of  $\theta $ such
that $ |\phi (X)| \leq  \tau $.

$L6.$ For every lower semi-continuous mapping $\theta :X \rightarrow
Y$ into a discrete space $Y$ there exists a set-valued selection $g:
X \rightarrow Y$ of  $\theta $ such that $\mid g (X)\mid \leq  \tau
$.

$L7.$ For every lower semi-continuous mapping $\theta :X \rightarrow
Y$ into a discrete space $Y$ there exists a single-valued selection
$g:X \rightarrow Y$ of  $\theta $ such that $\mid g (X)\mid \leq
\tau $.

$L8.$ Every open cover of $X$ has a subcover of cardinality $\leq
\tau. $ }

\textbf{Proof:}  Let $l(X) \leq \tau$. Then $k(X) \leq \tau ^+$ and
$\tau ^+$ is a regular cardinal. Theorem 2.3 completes the proof.
$\Box$

{\bf Theorem 2.5.} {\it Let $X$ be a regular space, $F$ be a compact
subset of $X$, $\tau $ be a   cardinal number and $k(Y) < \tau $ for
any closed subset $Y \subseteq  X\setminus F$ of $X$. Then
assertions $K1 - K8$ and $K12$ are equivalent. Moreover, if the
cardinal number $\tau $ is not sequential, then assertions $K1 -
K8$, $K11$ and $K12$ are equivalent.}

\textbf{Proof:}  Let $k(X) \leq \tau$, $F$ be a compact subset of
$X$, $\tau $ be a cardinal number and $k(Y) < \tau $ for any closed
subset $Y \subseteq  X\setminus F$ of $X$ and $\theta : X
\rightarrow Y$ be a lower semi-continuous closed-valued mapping into
a complete metric space $(Y, \rho )$.

Case 1. $\tau =\aleph_0$.

In this case the space $X$ is compact. Thus, from Theorem 0.1, it
follows that there exist a lower semi-continuous compact-valued
mapping $\varphi  :X \rightarrow Y$ and an upper semi-continuous
compact-valued mapping $\psi  :X \rightarrow Y$ such that $\varphi
(x) \subseteq \psi (x) \subseteq \theta (x)$ for any $x \in X.$ The
set $\psi (X)$ is compact and $\varphi (X) \subseteq \psi (X)$.
Implications $(K1\Rightarrow K9)$ and $(K1\Rightarrow K12)$ are
proved.

Case 2. $\tau $ is a regular cardinal number.

In this case Theorem 2.3 completes the proof.

Case 3. $\tau $ is an uncountable limit  cardinal.

Let $\tau ' = cf(\tau )$.

The subspace $F$ is compact. Thus, from the E.Michael's Theorem 0.1,
it follows that there exists an upper semi-continuous compact-valued
mapping $\psi  :F \rightarrow Y$ such that $ \psi (x) \subseteq
\theta (x)$ for any $x \in F.$ The set $\Phi = \psi (F)$ is compact.
There exists a sequence $\{H_n: n \in \mathbb N\}$ of open subsets
of $Y$ such that:

- $\Phi \subseteq H_{n+1} \subseteq cl_YH_{n+1} \subseteq H_n$ for
any $n \in \mathbb N;$

- for every open subset $V\supseteq \Phi $ of $Y$ there exists $n
\in \mathbb N$ such that $H_n \subseteq V.$

There exist a sequence $\gamma = \{\gamma _n =\{U_\alpha : \alpha
\in A_n\}: n \in \mathbb N\}$ of open covers of the space $X$, a
sequence $\xi  = \{\xi  _n =\{V_\alpha : \alpha \in A_n\}: n \in
\mathbb N\}$ of open families of the space $Y$, a sequence $\{U_n: n
\in \mathbb N\}$ of open subsets of $X$, a sequence $\pi = \{\pi _n:
A_{n+1} \rightarrow A_n: n \in \mathbb N\}$ of mappings  and a
sequence $\{\tau _n: n \in \mathbb N\}$ of cardinal numbers  such
that:

- $\cup \{U_\beta : \beta \in \pi _{n}^{-1}(\alpha )\}$ = $U_{\alpha
} \subseteq cl_X U_\alpha \subseteq \theta ^{-1}(V_\alpha )$ for any
$\alpha \in A_n$ and $n \in \mathbb N$;

- $\cup \{cl_{Y}V_\beta : \beta \in \pi _{n}^{-1}(\alpha )\}$
$\subseteq V_\alpha  $ and $diam(V_\alpha ) < 2^{-n}$ for any
$\alpha \in A_n$ and $n \in \mathbb N$;

- $|A_n| < \tau $ for any $n \in \mathbb N$;

- if $A_n' = \{\alpha \in A_n: F\cap cl_XU_\alpha =\emptyset \}$ and
$A_n'' = A_n\setminus A_n'$, then the set $A_n''$ is finite and $F
\subseteq U_n \subseteq cl_XU_n \subseteq \cup \{U_\alpha : \alpha
\in A_n''\}$;

- $\tau _{n} \leq \tau _{n+1} < \tau $ for any $n \in \mathbb N$;

- $cl_X U_n \subseteq \theta ^{-1}(H_n)$ and $|\{\alpha \in A_n:
U_\alpha \setminus  U_m \not= \emptyset \}|\leq \tau _m$ for all $n,
m \in \mathbb N;$

- $cl_XU_n \cap cl_XU_\alpha = \emptyset  $ for any $n \in \mathbb
N$ and $\alpha \in A_n'.$

Let $\eta = \{V: V$ {\it is open in $Y$ and} $diam (V) < 2^{-1}\}$.
There exists a finite subfamily $\{V_\beta  : \beta  \in B_1\}$ of
$\eta $ such that $\Phi \subseteq \cup \{V_\beta  : \beta  \in
B_1\}\subseteq H_1$. Let $W_1$ be an open subset of $Y$ and $\Phi
\subseteq W_1 \subseteq cl_Y W_1 \subseteq \cup \{V_\alpha : \alpha
\in B_1\}$.

Let $\gamma ' = \{U: U$ {\it is open in $X$ and $cl_{X}U \subseteq
\theta ^{-1}(V)$
 for some} $V \in \eta$  and $U \subseteq X\setminus U_1\}$ and
$\gamma '' = \{U: U$ {\it is open in $X$ and $cl_{X}U \subseteq
\theta ^{-1}(V_\beta )$
 for some} $\beta   \in B_1\}$.

Since $F$ is compact, there exist a finite family $\gamma _1'' =
\{U_\alpha : \alpha \in A_1''\}$ of $\gamma ''$ and an open subset
$U_1$ of $X$ such that $F\subseteq U_1 \subseteq cl_XU_1 \subseteq
\cup \{U_\alpha : \alpha \in A_1'' \}$ and $F\cap U_\alpha \not=
\emptyset $ for any $\alpha \in A_1''$. For every $\alpha \in A_1''$
fix $V_\alpha = V_\beta $ for some $\beta \in B_1$ such that
$clU_\alpha \subseteq \theta ^{-1}(V_\alpha ).$ Let $Y_1 =
X\setminus U_1.$ Since $k(Y) = \tau _1' < \tau ,$ there
 exists an open subfamily
$\gamma _1' =\{U_\alpha : \alpha \in A_1'\}$ of $\gamma '$ such that
$|A_1'| \leq  \tau _1,$
 $Y_1 \subseteq \cup \{U_\alpha : \alpha \in A_1'\}$ and
 $cl_XU_1 \cap (\cup \{cl_XU_\alpha : \alpha \in A_1'\} =\emptyset . $
For any $\alpha \in A_1'$ fix $V_\alpha \in \eta' $ such that $cl_X
U_\alpha \subseteq \theta ^{-1}(V_\alpha )$. Let $A_1 = A_1'\cup
A_1''$, $\gamma _1 =\{U_\alpha : \alpha \in A_1\}$ and $\eta  _1
=\{V_\alpha : \alpha \in A_1\}$.

Consider that the objects $\{\gamma _i, \xi _i, \pi _{i-1}, U_i,
\tau _i, : i \leq n\}$ are constructed.

We put $A_{im}= \{\alpha \in A_i: U_\alpha \cap U_m \not= \emptyset
\}$ for all $i,m \leq n.$

Fix $\alpha \in A_n$.

Let $\eta_{\alpha } = \{V: V$ {\it is open in $Y$, $cl_YV \subseteq
V_\alpha $
 and} $diam (V) < 2^{-n-1}\}$ and $\gamma_{\alpha }' = \{W: W$ {\it is open in $X$
and $cl_{X}W \subseteq \theta ^{-1}(V)$ for some} $V \in \eta_\alpha
\}$.

Assume that $\alpha \in A_n''$.

Since $F_\alpha = F\cap cl_XU_\alpha $ is a compact subset of $X$
there exists a finite subfamily $\gamma _{0\alpha} = \{W_\beta  :
\beta  \in A_{0\alpha }''\}$ of $\gamma_{\alpha }'$ such that
$F_\alpha \subseteq \cup \{W_\beta : \beta  \in A_{0\alpha }''\}$,
$F_\alpha \cap W_\beta \not= \emptyset $ for any $\beta  \in
A_{0\alpha }''$ and for any $\beta  \in A_{0\alpha}'' $ there exists
$V_\beta  \in \eta _\alpha $ such that $V_\beta \subseteq H_{n+1}$
and     $cl_X W_\beta  \subseteq \theta ^{-1}(V_\beta  )$. Now we
put $U_\beta = W_\beta \cap U_\alpha .$

Let $A_{n+1}'' = \cup \{A_{0\alpha }: \beta  \in A_n''\}$, $\gamma
_{n+1}'' = \{U_\beta  : \beta  \in A_{n+1}''\}$ and  $\eta  _{n+1}''
= \{V_\beta : \beta  \in A_{n+1}''\}$.

Let $\Phi _\alpha =cl_XU_\alpha \setminus \cup \{U_\beta : \beta \in
A_{0\alpha }''\}$ and $U_n' = U_n\setminus \cup \{\Phi _\alpha :
\alpha \in A_n''\}$. Then $U_n'$ is an open subset of $X$ and $F
\subseteq U_n' \cap \cup (\{U_\beta : \beta  \in A_{0\alpha }''\})$.

There exists an open subset $U_{n+1}$ of $X$ such that $U_{n+1}
\subseteq cl_XU_{n+1} \subseteq  U_n'\cap U_n \cap  (\cup \{U_\beta
: \beta  \in A_{0\alpha }''\})$.

Let $Y_i = X\setminus U_i$ for any $i \leq n+1.$ Then $\tau _i =
k(Y_i)$  for any $i \leq n+1.$

For any $\alpha \in A_n$  there exist the subfamilies $\gamma
_{i\alpha }' = \{W_\beta  : \beta  \in A_{in\alpha }'\}$, $i \leq
n+1,$ of $\gamma _\alpha '$ and the subfamilies $\eta  _{i\alpha }'
= \{V_\beta  : \beta \in A_{in\alpha }'\}$, $i \leq n+1,$ of $\gamma
_\alpha '$ such that:

- $|A_{in\alpha }'| < \tau _i$ for any $i \leq n+1;$

- $Y_i\cap cl_XU_\alpha \subseteq \cup \{W_\beta : \beta \in \cup
\{A_{jn\alpha }:j \leq i\}\}$ for any $i \leq n+1;$

- $Y_i\cap (\cup \{W_\beta : \beta \in \cup \{A_{jn\alpha }:i<j \leq
n+1\}\})= \emptyset $ for any $i < n+1.$

Now we put $A_{n\alpha } = \cup \{A_{in\alpha }: 0\leq i \leq
n+1\}$, $A_{n+1} = \cup \{A_{n\alpha }: \alpha \in A_n\}$, $U_\beta
= W_\beta \cap U_\alpha $, $\gamma _{n+1} = \{U_\beta :\beta  \in
A_{n+1}\}$, $\eta _{n+1} = \{V_\beta  :\beta  \in A_{n+1}\}$ and
$\pi _{n+1}^{-1}(\alpha ) = A_{n\alpha }$.

The objects $\{\gamma _n, \xi _n, \pi _{n}, U_n, \tau _n: n \in
\mathbb N\}$ are constructed.

Let $x \in X.$ Denote by $A(x)$ the set of all sequences $\alpha =
(\alpha _n: n \in \mathbb N)$ for which $\alpha _n \in A_n$ and $x
\in U_{\alpha _n}$ for any  $n \in \mathbb N.$ For any  $\alpha =
(\alpha _n: n \in \mathbb N) \in A(x)$ there exists a unique point
$y(\alpha ) \in Y$ such that $\{y(\alpha )\} = \cap \{V_{\alpha _n}:
n \in \mathbb N\}.$ It is obvious that $y(\alpha ) \in \theta (x).$
Let $\phi (x) = \{y(\alpha ): \alpha \in A(x)\}.$ Then $\phi $ is a
selection of $\theta $. By construction:

- $U_\alpha \subseteq \phi ^{-1}(V_\alpha )$ for all $\alpha \in
A_n$ and $n \in \mathbb N$;

- the mapping $\phi $ is lower semi-continuous;

- if $Z= \phi (X)$, then $\{H_\alpha = Z\cap V_\alpha : \alpha \in A
= \cup \{A_n:  n \in \mathbb N\}\}$ is an open base of the subspace
$Z$.

We affirm that $ k(cl_YZ) \leq \tau $.

Subcase 3.1. $\tau $ is not a sequential cardinal.

In this subcase $ m = sup \{|A_n|: n \in \mathbb N\} < \tau $ and
$w(Z) \leq |A| \leq m <\tau $. In this subcase we are proved the
implication $(K1 \rightarrow K11).$

Subcase 3.2. $\tau $ is a sequential cardinal.

Let $Z_n =\phi (Y_n)$ and $A_{nk} = \{\alpha \in A_n: Y_k\cap
U_\alpha \not=\emptyset \}$. Then $|Ank| < \tau _n$ for all $n, k
\in \mathbb N.$ Thus $w(Z_n)< \tau _n$.

Since $\phi (X)\setminus Z_n \subseteq H_n$, we have $k(cl_Y\phi (X)
\leq \tau .$

In this subcase we have proved implication $(K1 \rightarrow K2).$

Let $H = \cap \{U_n: n \in \mathbb N\}$, $\mu = \{\mu _n =\{H\cap
U_\alpha : \alpha \in A_n^{''}\}$ and $q \{q_n = \pi
_n|A_{n+1}^{''}: A_{n+1}^{''}\rightarrow A_{n}^{''} : n \in \mathbb
N\}$. By construction, we have
 $\cup \{W_\beta ; \beta \in q _{n}^{-1}(\alpha )\}$ =
$W_{\alpha } \subseteq cl_X W_\alpha \subseteq \theta ^{-1}(V_\alpha
)$ for any $\alpha \in A_n^{''}$ and $n \in \mathbb N$. Let $x \in
H.$ Denote by $B(x)$ the set of all sequences $\alpha = (\alpha _n:
n \in \mathbb N)$ for which $\alpha _n \in A_n^{''}$ and $x \in
cl_XW_{\alpha _n}$ for any  $n \in \mathbb N.$ For any  $\alpha =
(\alpha _n: n \in \mathbb N) \in B(x)$ there exists a unique point
$y(\alpha ) \in Y$ such that $\{y(\alpha )\} = \cap \{V_{\alpha _n}:
n \in \mathbb N\}.$ It is obvious that $y(\alpha ) \in \Phi \cap
\theta (x).$ Let $\mu _1 (x) = \{y(\alpha ): \alpha \in B(x)\}.$ The
mapping $\mu _1: H\rightarrow \Phi $ is compact-valued and upper
semi-continuous. Let $\mu (x) =\phi (x)$ for $x \in X\setminus H$
and $\mu (x) = \mu _1(x)$ for $x \in H$. Then $\mu  $ is a selection
of $\theta $. Fix a closed subset $Z \subseteq Y\setminus \Phi $ of
the space $Y$. Then $Z\cap \mu (X)\subseteq \phi (Y_n)$ for some $n
\in \mathbb N\}.$ Thus $w(Z\cap \mu (X) < \tau $. In this subcase we
have proved implication $(K1 \rightarrow K12),$ too.

Lemma 2.1 completes the proof of the theorem. $\Box$

The last theorem and Lemma 1.3 imply

{\bf Corollary 2.6.} {\it Let $X$ be a $\mu$-complete  space and
$\tau $ be a sequential  cardinal number.
 Then assertions $K1 - K8$ are equivalent.}

Theorem 2.5 is signigative for a sequential cardinal $\tau $. Every
compact subset of $X$ is paracompact in $X$. In fact we have

{\bf Theorem 2.7.} {\it Let $X$ be a regular space, $F$ be a
paracompact in $X$ subspace,
 $\tau $ be an infinite cardinal number, $k(F) \leq  \tau $,
$k(Y) < \tau $ for any closed subset $Y \subseteq  X\setminus F$ of
$X$.Then assertions $K1 - K8$  are equivalent.  Moreover, if the
cardinal number $\tau $ is not sequential, then assertions $K1 - K8$
and $K11$ are equivalent.}

\textbf{Proof:}  It is obvious that for any open in $X$ set
$U\supseteq F$ there exists an open subset $V$ of $X$ such that $F
\subseteq  U  \subseteq  cl_XU \subseteq  V$

Case 1. $\tau $ is a regular cardinal number.

In this case Theorem 2.3 completes the proof.

Case 2.  $\tau $ is a sequential cardinal number.

In this case Theorem 2.5 and Lemma 1.4 complete the proof.

Case 3. $\tau $ be a limit non-sequential  cardinal.

Let $\tau^* = cf(\tau )<\tau $. Obviously, $\tau^*$ is a regular
cardinal and $\tau ^* < \tau $.

There exist a sequence $\gamma = \{\gamma _n =\{U_\alpha : \alpha
\in A_n\}: n \in \mathbb N\}$ of open covers of the space $X$, a
sequence $\xi  = \{\xi  _n =\{V_\alpha : \alpha \in A_n\}: n \in
\mathbb N\}$ of open families of the space $Y$, a sequence $\{U_n: n
\in \mathbb N\}$ of open subsets of $X$ a sequence $\pi = \{\pi _n:
A_{n+1} \rightarrow A_n: n \in \mathbb N\}$ of mappings  and a
sequence $\{\tau _n: n \in \mathbb N\}$ of cardinal numbers  such
that:

- $\cup \{U_\beta ; \beta \in \pi _{n}^{-1}(\alpha )\}$ = $U_{\alpha
} \subseteq cl_X U_\alpha \subseteq \theta ^{-1}(V_\alpha )$ for
every $\alpha \in A_n$ and $n \in \mathbb N$;

- $\cup \{cl_{Y}V_\beta ; \beta \in \pi _{n}^{-1}(\alpha )\}$
$\subseteq V_\alpha  $ and $diam(V_\alpha ) < 2^{-n}$ for every
$\alpha \in A_n$ and $n \in \mathbb N$;

- $|A_n| < \tau _n \leq \tau _{n+1} < \tau$ for every $n \in \mathbb
N$;

- if $A_n' = \{\alpha \in A_n: F\cap cl_XU_\alpha =\emptyset \}$ and
$A_n'' = A_n\setminus A_n'$, then $\mid A_n''\mid <\tau^*$ and $F
\subseteq U_n \subseteq cl_XU_n \subseteq \cup \{U_\alpha : \alpha
\in A_n''\}$;

- the family $\gamma _n'' = \{U_\alpha :\alpha \in A_n''\}$ is
locally  finite in $X$ for every $n \in \mathbb N$;

- $cl_XU_n \cap cl_XU_\alpha = \emptyset  $ for every $n \in \mathbb
N$ and $\alpha \in A_n'.$

Let $\eta = \{V: V$ {\it is open in $Y$ and} $diam (V) < 2^{-1}\}$
and $\gamma ' = \{U: U$ {\it is open in $X$ and $cl_{X}U \subseteq
\theta ^{-1}(V))$ for some} $V \in \eta \}.$ There exist a locally
finite subfamily $\gamma _1'' = \{U_\alpha : \alpha \in A_1''\}$ of
$\gamma '$ such that $\mid A_1''\mid <\tau^* < k(F)$ and an open
subset $U_1$ of the  space $X$ such that $F\subseteq U_1 \subseteq
cl_XU_1 \subseteq \cup \{U_\alpha : \alpha \in A_1''\}$ and $F\cap
U_\alpha \not= \emptyset $ for every $\alpha \in A_1''$. For every
$\alpha \in A_1''$ fix $V_\alpha \in \eta $ such that $clU_\alpha
\subseteq \theta ^{-1}(V_\alpha ).$ Let $Y_1 = X\setminus U_1$ and
$\tau _1 = k(F)+\tau ^*$.  Since $k(Y) \leq \tau _1 < \tau ,$ there
 exists an open subfamily
$\gamma _1' =\{U_\alpha : \alpha \in A_1'\}$ of $\gamma '$ such that
$|A_1'| \leq  \tau _1,$
 $Y_1 \subseteq \cup \{U_\alpha : \alpha \in A_1'\}$ and
 $cl_XU_1 \cap (\cup \{cl_XU_\alpha : \alpha \in A_1'\}) =\emptyset . $
For every $\alpha \in A_1'$ fix $V_\alpha \in \eta' $ such that
$cl_X U_\alpha \subseteq \theta ^{-1}(V_\alpha )$. Let $A_1 =
A_1'\cup A_1''$, $\gamma _1 =\{U_\alpha : \alpha \in A_1\}$ and
$\eta  _1 =\{V_\alpha : \alpha \in A_1\}$.

The objects $\{\gamma _1, \xi _1, U_, \tau _1\}$ are constructed.

Consider that the objects $\{\gamma _i, \xi _i, \pi _{i-1}, U_i,
\tau _i, : i \leq n\}$ are constructed.

Fix $\alpha \in A_n$.

Let $\eta_{\alpha } = \{V: V$ {\it is open in $Y$, $cl_YV \subseteq
V_\alpha $ and} $diam (V) < 2^{-n-1}\}$ and $\gamma_{\alpha }^* =
\{W: W$ {\it is open in $X$ and $cl_{X}W \subseteq \theta ^{-1}(V)$
for some} $V \in \eta_\alpha   \}$.

Assume that $\alpha \in A_n''$.

Since $F_\alpha = F\cap cl_XU_\alpha $ is a closed subset of $X$,
then there exists a locally finite subfamily $\gamma _{\alpha}'' =
\{W_\beta  : \beta \in A_{\alpha }''\}$ of $\gamma_{\alpha }^*$,
where $\mid A_{\alpha }''\mid <\tau^*$ such that $F_\alpha \subseteq
\cup \{W_\beta : \beta \in A_{\alpha }''\}$, $F_\alpha \cap W_\beta
\not= \emptyset$ for every $\beta \in A_{\alpha }''\}$ and for every
$\beta  \in A_{\alpha}'' $ there exists $V_\beta  \in \eta _\alpha $
such that $cl_X W_\beta  \subseteq \theta ^{-1}(V_\beta  )$. We put
$U_\beta = U_\alpha \cap W_\beta $  for every $\beta  \in
A_{\alpha}''. $

Let $A_{n+1}'' = \cup \{A_{\alpha }'': \alpha \in A_n''\}$, $\gamma
_{n+1}'' = \{U_\alpha : \alpha \in A_{n+1}''\}$ and  $\eta  _{n+1}''
= \{V_\alpha : \alpha \in A_{n+1}''\}$.

The family $\gamma _{n+1}''$ is locally finite.

Let $\Phi _\alpha = cl_XU_\alpha \setminus \cup \{U_\beta :\beta \in
A_{\alpha }''$ and $U_n' = U_n\setminus \cup \{\Phi _\alpha : \alpha
\in A_n''\}$. Since the family $\gamma _n''$ is locally finite, the
set $U_n'$ is open in $X$ and $F \subseteq U_n' \subseteq \cup
\{U_\beta : \beta  \in A_{\alpha }''\}$.

There exists an open subset $U_{n+1}$ of $X$ such that $U_{n+1}
\subseteq cl_XU_{n+1} \subseteq \cup \{U_\beta : \beta  \in
A_{\alpha }''\}$.

Let $Y_{n+1} = X\setminus U_n$ and $\tau _{n+1} = k(Y_{n+1})+\tau
_n$.

For every $\alpha \in A_n$  there exist the subfamily $\gamma
_{\alpha }' = \{W_\beta  : \beta  \in A_{\alpha }'\}$
 of $\gamma _\alpha ^*$ and the subfamily $\eta  _{i\alpha }'
= \{V_\beta  : \beta \in A_{\alpha }'\}$ of $\gamma_\alpha '$ such
that:

- $|A_{\alpha }'| < \tau _{n+1};$

-$cl_XU_\alpha \setminus U_n \subseteq \cup \{W_\beta  : \beta \in
A_{\alpha }'\};$

- $cl_XW\beta \cap cl_XU_{n+1} =\emptyset$ for any $\beta \in
A_\alpha '.$

Now we put $A_{\alpha } = A_\alpha '\cup A_\alpha ''$,
 $A_{n+1} = \cup \{A_{\alpha }: \alpha \in A_n\}$,
$U_\beta = U_\alpha \cap U_\beta $ for any $\beta \in A_\alpha $,
$\gamma _{n+1} = \{U_\alpha :\alpha \in A_{n+1}\}$, $\eta _{n+1} =
\{V_\alpha :\alpha \in A_{n+1}\}$ and $\pi _{n+1}^{-1}(\alpha ) =
A_{n\alpha }$.

The objects $\{\gamma _n, \xi _n, \pi _{n}, U_n, \tau _n: n \in
\mathbb N\}$ are constructed.

Since $\tau $ is not sequential, we have $m = sup \{\tau _n: n \in
\mathbb N\} < \tau .$

Let $x \in X.$ Denote by $A(x)$ the set of all sequences $\alpha =
(\alpha _n: n \in \mathbb N)$ for which $\alpha _n \in A_n$ and $x
\in U_{\alpha _n}$ for every  $n \in \mathbb N.$ For every  $\alpha
= (\alpha _n: n \in \mathbb N) \in A(x)$ there exists a unique point
$y(\alpha ) \in Y$ such that $\{y(\alpha )\} = \cap \{V_{\alpha _n}:
n \in \mathbb N\}.$ It is obvious that $y(\alpha ) \in \theta (x).$
Let $\phi (x) = \{y(\alpha ): \alpha \in A(x)\}.$ Then $\phi $ is a
selection of $\theta $. By construction:

- $U_\alpha \subseteq \phi ^{-1}(V_\alpha ) $ for all $\alpha \in
A_n$ and $n \in \mathbb N$;

- the mapping $\phi $ is lower semi-continuous;

- if $Z= \phi (X)$, then $\{H_\alpha = Z\cap V_\alpha : \alpha \in A
= \cup \{A_n:  n \in \mathbb N\}\}$ is an open base of the subspace
$Z$ and $w(Z) \leq m$.

Thus we have proved the implication $(K1 \rightarrow K11).$

Lemma 2.1 completes the proof of the theorem. $\Box$

{\bf Remark 2.8.} {\it Let $X$ be a paracompact space and
$Y\subseteq X$. Then $l(cl_X Y)\leq l(Y)$ and $k(cl_X Y)\leq k(Y)$.}

Theorem 2.7, Corollary 2.6 and Lemma 1.4 yield

{\bf Corollary 2.9.} {\it Let $X$ be a paracompact and $\tau $ be an
infinite cardinal. Then the properties $K1 - K10$ are equivalent. }

One can observe that the Corollary 2.9 follows from Proposition 2.2,
Lemma 2.1 and Theorem 0.1, too.

{\bf Corollary 2.10.} {\it Let $X$ be a  space and $\tau$ be an
uncountable not sequential cardinal number. Then the following
assertions  are equivalent:

1. $X$ is a paracompact space and $k(X) \leq \tau$.

2. $X$ is a paracompact space and for every lower semi-continuous
closed-valued mapping $\theta :X\rightarrow Y$ into a complete
metrizable space $Y$ there exists a lower semi-continuous selection
$\phi:X\rightarrow Y$ of $\theta$ such that $w(\phi (X)) < \tau$.

3. $X$ is a paracompact space and for every lower semi-continuous
closed-valued mapping $\theta :X \rightarrow Y$ into a complete
metrizable space $Y$ there exists a single-valued selection $g :X
\rightarrow Y$ such that $w(g(X)) < \tau $.

4. $X$ is a paracompact space and for every lower semi-continuous
mapping $\theta :X \rightarrow Y$ into a discrete space $Y$ there
exists a single-valued selection $g :X \rightarrow Y$ such that
$\mid g (X) \mid < \tau $.

5. For every lower semi-continuous closed-valued mapping $\theta :X
\rightarrow Y$ into a complete metrizable space $Y$ there exist a
compact-valued lower semi-continuous mapping $\varphi  :X
\rightarrow Y$ and
 a compact-valued
upper semi-continuous mapping $\psi   :X \rightarrow Y$ such that
$w(\psi (X)) < \tau $ and $\varphi (x) \subseteq \psi (x) \subseteq
\theta (x)$ for any $x \in X.$

6 For every lower semi-continuous closed-valued mapping $\theta :X
\rightarrow Y$ into a complete metrizable space $Y$ there exists an
upper semi-continuous selection $\phi :X \rightarrow Y$ of  $\theta
$ such that $w(\phi (X)) < \tau $.}

{\bf Example 2.11.} Let $\tau $ be an uncountable limit cardinal
number and $m = cf(\tau )$. Fix a well ordered set $A$   and a
family of regular cardinal numbers $\{\tau_\alpha :\alpha\in A\}$
such that $\sup \{\tau_\alpha :\alpha \in A\}=\tau$ and
$\tau_\alpha< \tau _\beta  <\tau$ for all $\alpha , \beta \in A$ and
$\alpha < \beta $. For every $\alpha \in A$ fix a zero-dimensional
complete metric space $X_\alpha$ such that $w(X_\alpha)
=\tau_\alpha$. Let $X'$ be the discrete sum of the spaces
$\{X_\alpha:\alpha \in A\}$. Then $X'$ is a complete metrizable
space and $w(X')=\tau$. Thus $l(X')=\tau$ and $k(X')=\tau^+$. Fix a
point $b\not\in X'$. Put $X=\{b\}\cup X'$ with the topology
generated by the open bases $\{U\subseteq X': U$ {\it is open in}
$X'\} \bigcup \{X\setminus \bigcup\{X_\beta : \beta \leq \alpha \}:
\alpha \in A\}$. Then $X$ is a zero-dimensional paracompact space
and $\chi (X)=\chi (b,X)=cf(\tau)$. If $cf(\tau)=\aleph_0$, then $X$
is a complete metrizable space. If $Y \subseteq X'$ is a closed
subspace of $X$, then there exists $\alpha \in A$ such that $Y
\subseteq \cup \{X_\beta : \beta < \alpha \},$ $w(Y) < \tau _\alpha
$ and $k(Y) < \tau .$ Therefore $k(X) = \tau$.

Let $Z = X \times [0,1]$. Then $k(Z )= \tau $ and $k(Z,\tau ) =
\{b\}\times [0,1].$

Suppose that $\tau $ is not a sequential cardinal number, $\mathbb
N$ is a discrete space and $S = X\times \mathbb N$. Then $k(S) =
\tau $ and $k(S,\tau ) = \{b\}\times \mathbb N.$

Moreover, if $m = cf(\tau )$ is uncountable, $X_\tau $ is a complete
metrizable space, $w(X_\tau ) < m$ and $Z_\tau = X\times X_\tau $,
 then $k(Z_\tau ) = \tau $ and
$k(Z_\tau ,\tau ) = \{b\}\times X_\tau .$ $\Box$

\vskip 5pt

\section{On the geometry of paracompact spaces}

\vskip 5pt

Let $\Pi $ be the class of all paracompact spaces.

For every infinite cardinal number $\tau $ we denote by $\Pi (\tau
)$ the class $\{X \in \Pi : k(c\omega (X)) \leq \tau \}$. We put
$\Pi _l(\tau )$ = $\{X \in \Pi : l(c\omega (X)) \leq \tau \}.$

It is obvious that $\Pi (\tau ) \subseteq \Pi _l(\tau )\subseteq \Pi
(\tau ^+).$

We consider that $\Pi (n)=\{X \in \Pi : dim X = 0\}$ for any $n \in
\{0\}\cup \mathbb N$.

Our aim is to prove that the classes $\Pi (\tau )$ may be
characterized in terms of selections. The main results of the
section are the following two theorems.

{\bf Theorem 3.1.} {\it Let $X$ be a  space and $\tau$ be  an
uncountable non-sequential cardinal number. Then the following
assertions  are equivalent:

1. $X \in \Pi (\tau )$, i.e. $X$ is paracompact and $k(c\omega (X))
\leq  \tau$.

2. $X$ is a paracompact space and for every lower semi-continuous
closed-valued mapping $\theta :X\rightarrow Y$ into a complete
metrizable space $Y$ there exists a lower semi-continuous selection
$\phi:X\rightarrow Y$ of $\theta$ such that $w(\phi (c\omega (X))) <
\tau$.

3. $X$ is a paracompact space and for every lower semi-continuous
closed-valued mapping $\theta :X \rightarrow Y$ into a complete
metrizable space $Y$ there exists a single-valued selection $g :X
\rightarrow Y$ such that $w(g(c\omega (X))) < \tau $.

4. $X$ is a paracompact space and for every lower semi-continuous
mapping $\theta :X \rightarrow Y$ into a discrete space $Y$ there
exists a single-valued selection $g :X \rightarrow Y$ such that
$\mid g (c\omega (X)) \mid < \tau $.

5. For every lower semi-continuous closed-valued mapping $\theta :X
\rightarrow Y$ into a complete metrizable space $Y$ there exist a
compact-valued lower semi-continuous mapping $\varphi  :X
\rightarrow Y$ and
 a compact-valued
upper semi-continuous mapping $\psi   :X \rightarrow Y$ such that
$w(\psi (X)) < \tau $ and $\varphi (x) \subseteq \psi (x) \subseteq
\theta (x)$ for any $x \in c\omega (X).$

6. For every lower semi-continuous closed-valued mapping $\theta
:X\rightarrow Y$ into a complete metric space $Y$ there exist a
closed $G_\delta$-set $H$ of $X$ and an upper semi-continuous
compact-valued selection $\psi:X\rightarrow Y$ such that:

i)$c\omega (X) \subseteq H$ and $w(\psi (H))<\tau $;

ii) $\psi (x)$ is a one-point set of $Y$ for every $x\in X\setminus
H$;

iii) $cl_Y\psi (H) = cl_Y\psi (c\omega (X)).$

7. For every lower semi-continuous closed-valued mapping $\theta
:X\rightarrow Y$ into a complete metric space $Y$ there exists an
upper semi-continuous compact-valued selection $\psi: X\rightarrow
Y$ such that $k(\psi ^{\infty } (x))<\tau$ for every $x\in X$.

8. For every lower semi-continuous mapping $\theta :X\rightarrow Y$
into a discrete space $Y$ there exists an upper semi-continuous
selection $\psi:X\rightarrow Y$ such that $\mid \psi^{\infty} (x)
\mid <\tau $ for every $x\in X$.}

{\bf Theorem 3.2.} {\it Let $X$ be a  space and $\tau$ be  an
infinite cardinal number. Then the following assertions  are
equivalent:

1. $X \in \Pi (\tau )$, i.e. $X$ is paracompact and $k(c\omega (X))
\leq  \tau$.

2. $X$ is a paracompact space and for every lower semi-continuous
closed-valued mapping $\theta :X\rightarrow Y$ into a complete
metrizable space $Y$ there exists a lower semi-continuous selection
$\phi:X\rightarrow Y$ of $\theta$ such that $k(cl_Y\phi (c\omega
(X))) \leq  \tau$.

3. $X$ is a paracompact space and for every lower semi-continuous
closed-valued mapping $\theta :X \rightarrow Y$ into a complete
metrizable space $Y$ there exists a single-valued selection $g :X
\rightarrow Y$ such that $k(cl_Yg(c\omega (X))) \leq  \tau $.

4. $X$ is a paracompact space and for every lower semi-continuous
mapping $\theta :X \rightarrow Y$ into a discrete space $Y$ there
exists a single-valued selection $g :X \rightarrow Y$ such that
$\mid g (c\omega (X)) \mid < \tau $.

5. For every lower semi-continuous closed-valued mapping $\theta :X
\rightarrow Y$ into a complete metrizable space $Y$ there exist a
compact-valued lower semi-continuous mapping $\varphi  :X
\rightarrow Y$ and
 a compact-valued
upper semi-continuous mapping $\psi   :X \rightarrow Y$ such that
$k(cl_Y(\psi (c\omega (X))) \leq k(\psi (c\omega (X))) \leq  \tau $
and $\varphi (x) \subseteq \psi (x) \subseteq \theta (x)$ for any $x
\in c\omega (X).$

6. For every lower semi-continuous closed-valued mapping $\theta
:X\rightarrow Y$ into a complete metric space $Y$ there exist a
closed $G_\delta$-set $H$ of $X$ and an upper semi-continuous
compact-valued selection $\psi:X\rightarrow Y$ such that:

i)$c\omega (X) \subseteq H$ and $k(\psi (H))\leq \tau $;

ii) $\psi (x)$ is a one-point set of $Y$ for every $x\in X\setminus
H$;

iii) $cl_Y\psi (H) = cl_Y\psi (c\omega (X)).$

7. For every lower semi-continuous closed-valued mapping $\theta
:X\rightarrow Y$ into a complete metric space $Y$ there exists an
upper semi-continuous compact-valued selection $\psi: X\rightarrow
Y$ such that $k(\psi ^{n } (x))<\tau$ for every $x\in X$ and any $n
\in \mathbb N$.

8. For every lower semi-continuous mapping $\theta :X\rightarrow Y$
into a discrete space $Y$ there exists an upper semi-continuous
selection $\psi:X\rightarrow Y$ such that $\mid \psi^{n} (x) \mid
<\tau $ for every $x\in X$ and any $n \in \mathbb N$.}

\textbf{Proof of the Theorems:}
 Let  $X \in \Pi (\tau )$ and $\theta :X\rightarrow Y$ be a
lower semi-continuous closed-valued mapping into a complete metric
space $(Y,d)$. For every subset $L$ of $Y$ and every $n \in \mathbb
N$ we put $O(L,n) = \{y \in Y: d(y,L) = inf \{d(x,z): z \in L\} <
2^{-n}\}$. Obviously, $cl_YL$ = $\cap \{O(L,n): n \in \mathbb N$ and
$cl_YO(L,n+1) \subseteq  O(L,n)$ for any $n \in \mathbb N.$

By virtue of the Michael's Theorem 0.1, there exist a compact-valued
lower semi-continuous mapping $\varphi  :X \rightarrow Y$ and a
compact-valued upper semi-continuous mapping $\psi   :X \rightarrow
Y$ such that $\varphi (x) \subseteq \psi (x) \subseteq \theta (x)$
for any $x \in c\omega (X).$

>From Proposition 2.2 it follows that
$k(cl_Y(\psi (c\omega (X))) \leq k(\psi (c\omega (X))) \leq  \tau $
and $k(cl_Y(\varphi  (c\omega (X))) \leq k(cl_Y(\psi (c\omega (X))))
\leq  \tau $. Moreover, if $\tau $ is a not sequential cardinal
number, then $w(\varphi  (c\omega (X)) \leq w(\psi (c\omega (X))) <
\tau .$

Therefore, the assertions 2, 3, 4 and 5 of Theorems follow from the
assertion 1.

It will be affirmed that there exist a sequence
$\{\phi_n:X\rightarrow Y:\ n\in \mathbb{N}\}$ of lower
semi-continuous compact-valued mappings, a sequence
$\{\psi_n:X\rightarrow Y:n \in \mathbb{N}\}$ of upper
semi-continuous compact-valued mappings, a sequense $\{V_n: n
\mathbb N\}$ of open subsets of $Y$ and a sequense $\{H_n: n \mathbb
N\}$ of open-and-closed subsets of $X$   such that:

1) $\psi_{n+1}(x) \subseteq \phi_n (x) \subseteq \psi_n (x)
\subseteq \theta  (x)$ for every $x\in X$ and every $n \in
\mathbb{N}$;

2) $\phi_n (x)=\psi_n (x)$ is a one-point subset of $Y$ for every
$x\in X\setminus H_n$ and for every $n \in \mathbb{N}$;

3)$H_{n+1} \subseteq \{x\in X: \psi _n(x) \subseteq V_n$, $H_{n+1}
\subseteq H_n$ and $V_{n+1} = O(\psi _n(c\omega (X))$ for every $n
\in \mathbb{N}$;

Let $V_1 = O(\theta (c\omega (X))$ and $U_1 = \theta ^{-1}V_1$.
 From Lemma 0.2 it follows that there exists
an open-and-closed subset $H_1$ of $X$ such that $c\omega (X)
\subseteq H_{1}  \subseteq U_1$.

 Since $dim (X\setminus H_1)=0$ there
exists a single-valued continuous mapping $h_1:X\setminus
H_1\rightarrow Y$ such that $h_1 (x)\in \theta (x)$ for every $x\in
X\setminus H_1$. Since $H_1$ is a paracompat space, $V_1$ is a
complete metrizable space and $\theta _1: H_1 \rightarrow V_1$,
where $\theta _1(x) = V_1\cap \theta (x)$, is a lower semicontinuous
closed-valued in $V_1$ mapping, by virtue of Theorem 0.1, there
exist  a compact-valued lower semi-continuous mapping $\varphi  _1
:H_1 \rightarrow V_1$ and a compact-valued upper semi-continuous
mapping $\lambda _1  :H_1 \rightarrow V_1$ such that $\varphi _1 (x)
\subseteq \lambda _1 (x) \subseteq \theta _1(x)$ for any $x \in
H_1.$

Put $\psi_1 (x)=\phi _1 (x)= h_1(x)$ for $x \in X\setminus H_1$ and
$\psi_1 (x) = \lambda _1(x)$, $\phi _1(x) = \varphi _1(x)$ for $x
\in H_1.$

The objects $\phi_1$ and $\psi_1$  are constructed.

Suppose that $n > 1$ and the objects $\phi_{n-1}$,$ \psi_{n-1}$,
$H_{n-1}$ and $V_{n-1}$ had been constructed.

We put $F_n = cl_Y \psi _{n-1}(c\omega (X))$, $V_n= O(F_n, n)$ and
$U_n = \{x \in H_{n-1}: \psi _{n-1}(x) \subseteq V_n\}$.
 From Lemma 0.2 it follows that there exists
an open-and-closed subset $H_n$ of $X$ such that $c\omega (X)
\subseteq H_{n}  \subseteq U_n$.

 Since $dim (X\setminus H_{n})=0$ there
exists a single-valued continuous mapping $h_{n}:X\setminus
H_{n}\rightarrow Y$ such that $h_{n} (x)\in \phi _{n-1} (x)$ for
every $x\in X\setminus H_n$. By construction, we have $\phi _{n-1}
\subseteq \psi (x) \subseteq V_n $ for any $x \in H_n.$
 Since $H_n$ is a paracompat space,
$V_n$ is a complete metrizable space and $\theta _n: H_n \rightarrow
V_n$, where $\theta _n(x) = V_n\cap \phi _{n-1} (x)$, is a lower
semicontinuous closed-valued in $V_n$ mapping, by virtue of Theorem
0.1, there exist  a compact-valued lower semi-continuous mapping
$\varphi  _n  :H_n \rightarrow V_n$ and a compact-valued upper
semi-continuous mapping $\lambda _n  :H_n \rightarrow V_n$ such that
$\varphi _n (x) \subseteq \lambda _n (x) \subseteq \theta _n(x)$ for
any $x \in H_n.$

Put $\psi_n (x)=\phi _n (x)= h_n(x)$ for $x \in X\setminus H_n$ and
$\psi_n (x) = \lambda _n(x)$, $\phi _n(x) = \varphi _n(x)$ for $x
\in H_n.$ The objects $\phi_n$ and $\psi_n$  are constructed.

Now we put $\lambda (x) = \cap \{\psi _n(x): n \in \mathbb N\}$ for
any $x \in X $  and $H = \cap \{H_n: n \in \mathbb N\}$.

Sinse $\lambda ^{-1}(\Phi ) = \cap \{\psi _n^{-1}(\Phi ): n \in
\mathbb N\}$ for any closed subset $\Phi $ of $Y$, the mapping
$\lambda $ is compact-valued and upper semi-continuous. By
construction,

i) $c\omega (X) \subseteq H$ and $k(\lambda  (H))\leq \tau $;

ii) $\lambda  (x)$ is a one-point set of $Y$ for every $x\in
X\setminus H$;

iii) $cl_Y\lambda  (H) = cl_Y \lambda  (c\omega (X));$

iv) $\lambda (\lambda ^{-1}(A) \subseteq A\cup \lambda (H)$ for
every subset $A$ of $Y$.

Therefore,  the assertions 6, 7 and 8 of Theorems follow from the
assertion 1.

$(8\Rightarrow 1)$ Let $\gamma  =\{U_\alpha:\ \alpha \in A\}$ be an
open cover of $X$. On $A$ introduce the discrete topology and put
$\theta (x)= \{\alpha \in A:\ x \in U_\alpha \}$ for $x\in X$. Since
$\theta^{-1} (H)=\bigcup \{U_\alpha:\ \alpha \in H\}$ for every
subset $H$ of $A$, the mapping $\theta :X\rightarrow A$ is lower
semi-continuous. Let $\psi :X\rightarrow A$ be an upper
semi-continuous selection of $\theta $ with $\mid \psi^2 (x)\mid
<\tau$ for every $x\in X$. Then $\xi  =\{\Psi_\alpha = \psi^{-1}
(\alpha ):\ \alpha \in A\}$ is a closed closure-preserving
$\tau^-$-star shrinking of the cover $\xi$.By virtue of Proposition
1.5, the assertion 1 follows from the assertion 8. $\Box$

{\bf Corollary 3.3.} {\it For a topological space $X$ the following
assertions  are equivalent:

1) $X$ is paracompact and $c\omega (X)$ is compact.

2) $X$ is strongly paracompact and $c\omega (X)$ is compact.

3) For every lower semi-continuous closed-valued mapping $\theta
:X\rightarrow Y$ into a complete metric space $Y$ there exist an
upper semi-continuous compact-valued selection $\psi :X\rightarrow
Y$ and a closed $G_\delta$-subset $H$ of $X$ such that $c\omega (X)
\subseteq H$, $cl_Y (\psi (H))$ is compact and $\psi (x)$ is a
one-point set for every $x\in X\setminus H$.

4) For every lower semi-continuous closed-valued mapping $\theta
:X\rightarrow Y$ into a complete metric space $Y$ there exists an
upper semi-continuous selection $\psi :X\rightarrow Y$ such that
$cl_Y \psi^\infty (x)$ is compact for every $x\in X$.

5) For every lower semi-continuous mapping $\theta :X\rightarrow Y$
into a discrete space $Y$ there exists an upper semi-continuous
selection $\psi :X\rightarrow Y$ such that the set $\psi^\infty (x)$
is finite for every $x\in X$.

6) For every open cover of $X$ there exists an open star-finite
shrinking.}

\textbf{Proof:} For the implication $(1\Rightarrow 2)$ see
Proposition 4, \cite{ChMN}.

For the implications $(1\Leftrightarrow 6)$ see Proposition 5,
\cite{ChMN}. $\Box$

\vskip 5pt

{\bf Corollary 3.4.} {\it For a  space and an infinite cardinal
number $\tau$ the following assertions  are equivalent:

1) $X$ is paracompact and $l(c\omega (X))\leq \tau $.

2) For every lower semi-continuous closed-valued mapping $\theta
:X\rightarrow Y$ into a complete metric space $Y$ there exist an
upper semi-continuous compact-valued selection $\psi :X\rightarrow
Y$ and a closed $G_\delta$-subset $H$ of $X$ such that $c\omega (X)
\subseteq H$ and $w( \psi (H))\leq \tau$; $\psi (x)$ is a one-point
set for every $x\in X\setminus H$.

3) For every lower semi-continuous closed-valued mapping $\theta
:X\rightarrow Y$ into a complete metric space $Y$ there exists an
upper semi-continuous compact-valued selection $\psi :X\rightarrow
Y$ such that $w(\psi^\infty (x))\leq \tau$ for every $x\in X$.

4) For every lower semi-continuous mapping $\theta :X\rightarrow Y$
into a discrete space $Y$ there exists an upper semi-continuous
selection $\psi :X\rightarrow Y$ such that $\mid \psi^\infty (x)\mid
\leq \tau$ for every $x\in X$. }

{\bf Corollary 3.5.} {\it For a topological space $X$ the following
assertions  are equivalent:

1) $X$ is paracompact and $c\omega (X)$ is Lindel\"{o}f.

2) $X$ is strongly paracompact and $c\omega (X)$ is Lindel\"{o}f.

3) For every lower semi-continuous closed-valued mapping $\theta
:X\rightarrow Y$ into a complete metric space $Y$ there exist an
upper semi-continuous compact-valued selection $\psi :X\rightarrow
Y$ and a closed $G_\delta$-subset $H$ of $X$ such that $c\omega (X)
\subseteq H$, $\psi (H)$ is separable and $\psi (x)$ is a one-point
set for every $x\in X\setminus H$.

4) For every lower semi-continuous closed-valued mapping $\theta
:X\rightarrow Y$ into a complete metric space $Y$ there exists an
upper semi-continuous compact-valued selection $\psi :X\rightarrow
Y$ such that $\psi^\infty (x)$ is separable for every $x\in X$.

5) For every lower semi-continuous mapping $\theta :X\rightarrow Y$
into a discrete space $Y$ there exists an upper semi-continuous
selection $\psi :X\rightarrow Y$ such that the set $\psi^\infty (x)$
is countable for every $x\in X$.

6) For every open cover of $X$ there exists an open star-countable
shrinking.}

{\bf Example 3.6.} Let $A$ be an uncountable set and $X_\alpha $ be
a non-empty compact space for every $\alpha \in A$. Let $X=\bigoplus
\{X_\alpha:\ \alpha \in A\}$ be the discrete sum of the space
$\{X_\alpha:\ \alpha \in A\}$. Let $B=\{\alpha \in A:\ dim X_\alpha
\neq 0\}$. Then $c\omega (X)$ is compact if and only if the set $B$
is finite. If the set $B$ is infinite then $l(c\omega (X))=\mid
B\mid$ and $k(c\omega (X))=\mid B\mid^+$.$\Box$

{\bf Example 3.7.} Let $\tau $ be an uncountable non-sequential
cardinal number. Fix an infinite set $A_m$ for every cardinal number
$m<\tau$ assuming that $A_m \cap A_n = \emptyset$ for $m\neq n$. Put
$A= \bigcup \{A_m:\ m<\tau\}$. Let $\{X_\alpha :\ \alpha \in A\}$ be
a family of non-empty compact spaces assuming that $X_\alpha \cap
X_\beta = \emptyset $ for $\alpha \neq \beta$. Put $B_m=\{\alpha \in
A_m:\ dim X_\alpha \neq 0\}$ and $1\leq \mid B_m \mid \leq m$ for
every $m<\tau$. Fix a point $b\not\in \bigcup \{X_\alpha:\ \alpha
\in A\}$. Let $X=\{b\}\cup (\bigcup \{X_\alpha :\ \alpha \in A\})$.
Suppose that $X_\alpha $ is an open subset of $X$ and
$\{H_m=\{b\}\cup (\bigcup \{X_\alpha :\ \alpha \in A_n,\ n\leq
m\}):\ m<\tau\}$ is a base of $X$ at $b$. If $Z=\{b\}\cup (\bigcup
\{X_\alpha :\ \alpha \in B_m,\ m<\tau\})$, then $c\omega
(X)\subseteq Z$ and $k(c\omega (X))\leq k(Z)=l(Z)=\tau$. $\Box$

{\bf Example 3.8.} Let $\tau$ be a regular uncountable cardinal
number, $A$ be an infinite set, $\tau<\mid A \mid $, $\{X_\alpha :\
\alpha \in A\}$ be a family of non-empty compact spaces, $X_\alpha
\cap X_\beta = \emptyset$ for $\alpha \neq \beta ,\ B=\{\alpha \in
A:\ dim X_\alpha \neq 0\}, \ \tau = \mid B\mid $ and $b\not\in
\bigcup \{X_\alpha:\ \alpha \in A\}$. Let $X=\{b\}\cup (\bigcup
\{X_\alpha :\ \alpha \in A\})$. Suppose that $X_\alpha $ is an open
subset of $X$ and $\{U_H=X \setminus \bigcup \{X_\alpha :\ \alpha
\in H\}:\ H \subseteq A, \ \mid H \mid <\tau \}$ is a base of $X$ at
$b$. If $Z=\{b\}\cup (\bigcup \{X_\alpha :\ \alpha \in B\})$, then
$c\omega (X)\subseteq Z$ and $k(c\omega (X))\leq k(Z)=l(Z)=\tau$.
$\Box$

{\bf Example 3.9.} Let $\tau$ be a regular uncountable limit
cardinal number and $2^m < \tau $ for any $m <\tau .$ Let
$\{m_\alpha : \alpha \in A\}$ be a family of infinite cardinal
numbers such that $|A| = \tau $, the set $A$ is well ordered and
$m_\alpha < m_\beta  $, $|\{\mu \in A: \mu \leq \alpha \}| < \tau $
provided $\alpha , \beta  \in A$ and $\alpha <\beta .$ For any
$\alpha \in A$ fix a discrete space of the cardinality $m_\alpha .$
Let $X = \Pi \{X_\alpha : \alpha \in A\}$. If $x = (x_\alpha :\alpha
\in A) \in X$ and $\beta \in A$, then $O(\beta ,x) = \{y = (y_\alpha
:\alpha \in A) \in X: y_\alpha = x_\alpha$ {\it for any} $\alpha
\leq \beta \}$. The family $\{O(\beta ,x): \beta \in A, x \in X\}$
form the open base of the space $X$. The space $X$ is paracompact
and $w(X) = l(X) = \tau $. It is obvious that $c(X,\tau ) = X$, we
have $k(X) = \tau ^+.$ If $\alpha \in A$, then $\gamma _\alpha  =
\{O(\alpha ,x): x \in X\}$ is open discrete cover of $X$ and $|
\gamma _\alpha | = 2^{m_\alpha } < \tau $. $\Box$

\vskip 5pt

\section {On the class $\Pi (0)$ of spaces}

\vskip 5pt

In the present section the class of all paracompact spaces $X$ such
that $dim X =0$ is studied.

{\bf Definition 4.1} {\it A set-valued mapping  $\psi
:X\longrightarrow Y$ {\it is called virtual single-valued if
}$\psi^\infty (x)=\psi (x)$  for every } $x\in X$.

{\bf Remark 4.2}
 {\it It is obvious that for a set-valued mapping $\theta :X\longrightarrow Y$
the following conditions are equivalent:

1. $\psi  $ is a virtual single-valued mapping;

2. $\psi^2 (x)=\psi (x)$ for every $x\in X\ $;

3. $\psi^n (x)=\psi (x)$ for every $x\in X\ $  and some $n\geq 2$;

4.  $\psi (x)= \psi (y)$ provided $x, y\in X$ and $\psi (x)\cap \psi
(y)\neq \emptyset$.

5.  $\psi ^{-1}(y)= \psi ^{-1}(z)$ provided $y, z\in Y$ and $\psi
^{-1}(y)\cap \psi ^{-1}(z)\neq \emptyset$}.

Note that, if $f:X\longrightarrow Y$ is a single-valued mapping onto
a space $Y$, then $f^{-1}$ and $f$  are virtual single-valued
mappings.

Denote with $D=\{0,1\}$ the two-point discrete space.

{\bf Theorem 4.3} {\it For a space $X$, the following assertions
are equivalent:

1. $X$ is normal and $dimX=0$;

2. For every lower semi-continuous mapping $\theta :X\longrightarrow
D$ there exists a virtual single-valued lower semi-continuous
selection;

3. For every lower semi-continuous mapping $\theta :X\longrightarrow
D$ there exists a virtual single-valued upper semi-continuous
selection;

4. For every lower semi-continuous mapping $\theta :X\longrightarrow
D$  there exists a single-valued continuous selection. }

\textbf{ Proof: } Implications $(1\Leftrightarrow 4)$ is a well
known fact. Implications $(4\Rightarrow 2)$ and $(4\Rightarrow 3)$
are obvious as every single-valued continuous selection is virtual
single-valued.

$(2\Rightarrow 1)$ and $(3\Rightarrow 1)$ Let $F_1$ and $F_2$ be two
disjoint closed subsets of $X$. Put $\theta (x)= \{0\}$ for $x\in
F_1$, $\theta (x)= \{1\}$ for $x\in F_2$ and $\theta (x)= \{0,1\}$
for $x\in X\setminus (F_1 \cup F_2).$ The mapping $\theta
:X\longrightarrow D$ is lower semi-continuous. Suppose that $\lambda
:X\longrightarrow D$ is a virtual single-valued selection of
$\theta$. Put $H_1=\lambda^{-1}(0)$ and $H_2=\lambda^{-1}(2)$. Then
$F_1\subseteq H_1$ and $F_2\subseteq H_2$., $X=H_1\cup H_2 $ and
$H_1\cap H_2 =\emptyset$. If $\lambda $ is lower semi-continuous (or
upper semi-continuous)) the sets $H_1, H_2$ are open (closed).
$\Box$

Let $\tau $ be an infinite cardinal number. A topological space $X$
is called \emph{$\tau $-paracompact} if $X$ is normal and every open
cover of $X$ of the cardinality $\leq \tau $ has a locally finite
open refinement.

{\bf Theorem 4.4}
 {\it For a space $X$ and an infinite cardinal number $\tau $
the following assertions are equivalent:

 1. $X$ is a $\tau $-paracompact space and $dim X= 0.$

 2.For every lower semi-continuous mapping
$\theta :X\longrightarrow Y$ into a complete metrizable space $Y$ of
the weight $\leq \tau $ there exists a virtual single-valued lower
semi-continuous selection;

3. For every lower semi-continuous mapping $\theta :X\longrightarrow
Y$ into a complete metrizable space $Y$ of the weight $\leq \tau $
there exists a virtual single-valued upper semi-continuous
selection;

4. For every lower semi-continuous mapping $\theta :X\longrightarrow
Y$ into a complete metrizable space $Y$ of the weight $\leq \tau $
there exists a  single-valued continuous selection;

5. For every lower semi-continuous mapping $\theta :X\longrightarrow
Y$  into a discrete space $Y$ of the cardinality $\leq \tau $ there
exists a single-valued continuous selection. }

\textbf{ Proof: } Let $\gamma = \{U_\alpha : \alpha \in A\}$ be an
open cover of $X$ and $|A|\leq  \tau $. Consider that $A$ is a
wellordered discrete space and $\theta (x) = \{\alpha \in A: x \in
U_\alpha \}$ for any $x \in X.$ Then $\theta $ is a lower
semi-continuous mapping. Suppose that $\psi : X \rightarrow Y$ is a
a virtual single-valued lower or upper semi-continuous selection of
$\theta $. For any $x \in X$ we denote by $f(x)$ the first element
of the set $\psi (x)$. Then $f: X \rightarrow Y$ is a single-valued
continuous selection of the mappings $\theta $ and $\psi $.
Therefore
 $H_\alpha =  f^{-1}(\alpha ): \alpha \in A\}$
is a discrete refinement of $\gamma $. The implications
$(2\Rightarrow 1)$, $(2\Rightarrow 4)$, $(3\Rightarrow 1)$,
$(3\Rightarrow 4)$ and $(5\Rightarrow 1)$ are proved. The
implications $(4\Rightarrow 5)$, $(4\Rightarrow 2)$ and
$(4\Rightarrow 3)$ are obvious. The implication $(1\Rightarrow 4)$
is wellknown (see [1, 2]). $\Box$

{\bf Corollary 4.5}
 {\it For a space $X$  the following assertions are equivalent:

 1. $X$ is a paracompact space and $dim X= 0.$

 2. For every lower semi-continuous mapping
$\theta :X\longrightarrow Y$ into a complete metrizable space $Y$
there exists a virtual single-valued lower semi-continuous
selection;

3. For every lower semi-continuous mapping $\theta :X\longrightarrow
Y$ into a complete metrizable space $Y$ there exists a virtual
single-valued upper semi-continuous selection;

4. For every lower semi-continuous mapping $\theta :X\longrightarrow
Y$ into a complete metrizable space $Y$ there exists a
single-valued continuous selection;

5. For every lower semi-continuous mapping $\theta :X\longrightarrow
Y$  into a discrete space $Y$ there exists a single-valued
continuous selection. }

{\bf Remark 4.6}  {\it Let $Y$ be a topological space. Then:

1. If the space $Y$ is discrete, then every lower semi-continuous
virtual single-valued mapping or every upper semi-continuous virtual
single-valued mapping $\theta :X\longrightarrow Y$  into the space
$Y$  is continuous.

1. If the space $Y$ is not discrete, then there exist a paracompact
space $X$ and a virtual single-valued  mapping $\theta
:X\longrightarrow Y$ such that;

- $\theta $ is upper semi-continuous and not continuous;

- $X$ has a unique not isolated point.

3. If  $Y$ has an open non-discrete subspace $U$ and $|U| \leq
|Y\setminus U|$, then there exist a paracompact space $X$ and a
virtual single-valued  mapping $\theta :X\longrightarrow Y$ such
that;

- $\theta $ is lower semi-continuous and not continuous;

- $X$ has a unique not isolated point. }

{\bf Remark 4.7}  {\it Let $\gamma = \{Hy: y \in Y\}$ be a cover of
a space $X$, $Y$ be a discrete space and $\theta _\gamma (x) = \{y
\in Y: x \in Hy\}$. Then:

- the mapping $\theta _\gamma $ is lower semi-continuous if and only
if $\gamma $ is an open cover;

- the mapping $\theta _\gamma $ is upper semi-continuous if and only
if $\gamma $ is a closed and conservative cover;

- the mapping  $\varphi : X \rightarrow Y$ is a selection of  the
mapping $\theta _\gamma $
 if and only if $\{ Vy =\varphi ^{-1}(y): y \in Y\}$ is a shrinking of $\gamma $.}

Therefore, the study of the problem of the selections for the
mappings into discrete spaces is an  essencial case of the this
problem.

\end{document}